\documentclass[12pt,a4paper]{article}

\usepackage[margin=2.5cm]{geometry}
\usepackage{listings}
\usepackage{color}
\usepackage{textcomp}
\usepackage{tikz}
\usepackage{float}
\usepackage{amsthm}
\usepackage{amsmath}
\usepackage{amssymb}
\usepackage{mathrsfs}
\usepackage{pythonhighlight}
\usepackage{dsfont}
\usepackage{enumitem}
\usepackage{graphicx}      
\usepackage{subcaption}    
\usepackage[authoryear,round]{natbib} 
\usepackage{authblk}       
\usepackage{hyperref}      

\theoremstyle{plain}
\newtheorem{theorem}{Theorem}[section]
\newtheorem{lemma}[theorem]{Lemma}

\newtheorem{proposition}[theorem]{Proposition}
\newtheorem{corollary}[theorem]{Corollary}
\theoremstyle{remark}
\newtheorem{definition}[theorem]{Definition}
\newtheorem{remark}[theorem]{Remark}

\newlength{\starsize}
\newlength{\starspread}
\tikzset{starsize/.code={\setlength{\starsize}{#1}},
	starspread/.code={\setlength{\starspread}{#1}}}
\tikzset{starsize=1mm,
	starspread=3mm}
\pgfdeclarepatternformonly[\starspread,\starsize]
{custom fivepointed stars}
{\pgfpointorigin}
{\pgfqpoint{\starspread}{\starspread}}
{\pgfqpoint{\starspread}{\starspread}}
{
	\pgftransformshift{\pgfqpoint{\starsize}{\starsize}}
	\pgfpathmoveto{\pgfqpointpolar{18}{\starsize}}
	\pgfpathlineto{\pgfqpointpolar{162}{\starsize}}
	\pgfpathlineto{\pgfqpointpolar{306}{\starsize}}
	\pgfpathlineto{\pgfqpointpolar{90}{\starsize}}
	\pgfpathlineto{\pgfqpointpolar{234}{\starsize}}
	\pgfpathclose%
	\pgfusepath{fill}
}
\usetikzlibrary{backgrounds}
\usetikzlibrary{patterns,fadings}
\usetikzlibrary{arrows,decorations.pathmorphing}
\usetikzlibrary{calc}
\definecolor{light-gray}{gray}{0.95}
\newlength{\hatchspread}
\newlength{\hatchthickness}
\newlength{\hatchshift}
\newcommand{\hatchcolor}{}
\tikzset{hatchspread/.code={\setlength{\hatchspread}{#1}},
	hatchthickness/.code={\setlength{\hatchthickness}{#1}},
	hatchshift/.code={\setlength{\hatchshift}{#1}},
	hatchcolor/.code={\renewcommand{\hatchcolor}{#1}}}
\tikzset{hatchspread=3pt,
	hatchthickness=0.4pt,
	hatchshift=0pt,
	hatchcolor=black}
\pgfdeclarepatternformonly[\hatchspread,\hatchthickness,\hatchshift,\hatchcolor]
{custom north west lines}
{\pgfqpoint{\dimexpr-2\hatchthickness}{\dimexpr-2\hatchthickness}}
{\pgfqpoint{\dimexpr\hatchspread+2\hatchthickness}{\dimexpr\hatchspread+2\hatchthickness}}
{\pgfqpoint{\dimexpr\hatchspread}{\dimexpr\hatchspread}}
{
	\pgfsetlinewidth{\hatchthickness}
	\pgfpathmoveto{\pgfqpoint{0pt}{\dimexpr\hatchspread+\hatchshift}}
	\pgfpathlineto{\pgfqpoint{\dimexpr\hatchspread+0.15pt+\hatchshift}{-0.15pt}}
	\ifdim \hatchshift > 0pt
	\pgfpathmoveto{\pgfqpoint{0pt}{\hatchshift}}
	\pgfpathlineto{\pgfqpoint{\dimexpr0.15pt+\hatchshift}{-0.15pt}}
	\fi
	\pgfsetstrokecolor{\hatchcolor}
	\pgfusepath{stroke}
}
\pgfdeclarepatternformonly[\hatchspread,\hatchthickness,\hatchshift,\hatchcolor]
{custom north east lines}
{\pgfqpoint{\dimexpr-2\hatchthickness}{\dimexpr-2\hatchthickness}}
{\pgfqpoint{\dimexpr\hatchspread+2\hatchthickness}{\dimexpr\hatchspread+2\hatchthickness}}
{\pgfqpoint{\dimexpr\hatchspread}{\dimexpr\hatchspread}}
{
	\pgfsetlinewidth{\hatchthickness}
	\pgfpathmoveto{\pgfqpoint{\dimexpr\hatchshift-0.15pt}{-0.15pt}}
	\pgfpathlineto{\pgfqpoint{\dimexpr\hatchspread+0.15pt}{\dimexpr\hatchspread-\hatchshift+0.15pt}}
	\ifdim \hatchshift > 0pt
	\pgfpathmoveto{\pgfqpoint{-0.15pt}{\dimexpr\hatchspread-\hatchshift-0.15pt}}
	\pgfpathlineto{\pgfqpoint{\dimexpr\hatchshift+0.15pt}{\dimexpr\hatchspread+0.15pt}}
	\fi
	\pgfsetstrokecolor{\hatchcolor}
	\pgfusepath{stroke}
}

\DeclareMathOperator{\spann}{span}
\newcommand{\dd}{\text{\rm d}}
\newcommand{\msf}[1]{{\mathsf #1}}
\newcommand{\mc}[1]{{\mathcal #1}}
\newcommand{\mf}[1]{{\mathfrak #1}}

\newcommand{\bb}[1]{{\mathbb #1}}

\newcommand{\norm}[1]{\lVert#1\rVert}
\newcommand{\p}{\partial}

\def\centerarc[#1](#2)(#3:#4:#5){\draw[#1] ($(#2)+({#5*cos(#3)},{#5*sin(#3)})$) arc (#3:#4:#5);}
\newcommand{\pfrac}[2]{\genfrac{}{}{}{1}{#1}{#2}}
\newcommand{\mud}[1]{\left\lvert#1\right\rvert}
\newcommand{\one}{\mathds{1}}
\newcommand{\pn}[1]{\big(\pfrac{#1}{n}\big)}

\begin{document}
	
	\title{A functional central limit theorem for the general Brownian motion on the half-line}

	\author[1]{Dirk Erhard\thanks{dirk.erhard@ufba.br}}
	\author[1]{Tertuliano Franco\thanks{tertu@ufba.br}}
	\author[2]{Milton Jara\thanks{mjara@impa.br}}
	\author[1]{Eduardo Pimenta\thanks{eduardo.sampaio@ufba.br}}

	\affil[1]{UFBA, Instituto de Matem\'atica e Estat\'istica, Campus de Ondina, Av. Milton Santos, S/N. CEP 40170-110, Salvador, Brazil}
	\affil[2]{IMPA, Estrada Dona Castorina, 110, Jardim Bot\^anico, CEP 22460-320, Rio de Janeiro, Brazil}

	\date{}

	\maketitle

	\begin{abstract}
		In this work, we establish a Trotter-Kato type theorem. More precisely, we characterize the convergence in distribution of Feller processes by examining the convergence of their generators. The main novelty lies in providing quantitative estimates in the vague topology at any fixed time. As important applications, we deduce functional central limit theorems for random walks on the positive integers with boundary conditions, which converge to Brownian motions on the positive half-line with boundary conditions at zero.
	\end{abstract}




	\section{Introduction}\label{s1}
	
	The subject of functional central limit theorems (functional CLTs) originated from the now standard Donsker Theorem and Invariance Principles for Brownian motion. Since then, a substantial body of literature has emerged, focusing on invariance principles for various types of random walks (including those on random media) and their convergence to standard Brownian motion.
	However, the convergence to Brownian-type processes, such as the skew Brownian motion, sticky Brownian motion, elastic Brownian motion, and others, has received much less attention up to the present day.
	
	As some rare examples, \cite{Harrison_Lemoine} showed the convergence of a specific queue to sticky Brownian motion. \cite{Amir} established the convergence of rescaled discrete-time random walks with deterministic waiting times, also to sticky Brownian motion. More recently, \cite{EFS2020} proved a functional central limit for the slow bond random walk. The limit of this process is given by the \textit{snapping out Brownian motion}, a Brownian-type process created  by \cite{Lejay}. Somewhat related to the subject,  \cite{cookie} have shown the convergence of the one-dimensional cookie random walk (a walk that takes a decision based on the local time of the present position) towards what they called a \textit{Brownian motion perturbed at extrema}, which is a stochastic process $W(t)$ solving a functional equation relating $W(t)$, its maxima and minima, and a standard Brownian motion.
	
	In this work, we present a criterion to ensure a functional central limit theorem for Feller processes, based on the rate of convergence of their generators. This criterion comes with a corresponding Berry-Esseen type estimate. Convergence in distribution via the convergence of generators is not a new topic. A classical result in the book \citet[Theorem 6.1 on page 28]{EK} can be paraphrased as follows: Consider a sequence $(T_n)$ of strongly continuous contraction semigroups, where each $T_n$ is defined on a Banach space $\text{X}_n$ and has generator $L_n$. Denote by $\pi_n$ the projection of another Banach space $\text{X}$ onto $\text{X}_n$ where a strongly continuous semigroup $T$ with generator $L$ is defined. Then $T_n(t)\pi_n f$ converges to $T(t)f$ for each $t$ if and only if for every $f$ in a core for $L$ there exists $f_n$ in a core for $L_n$ such that $f_n\to f$ and $L_n f_n\to Lf$. Our main result complements the above result and shows that a rate of convergence in the convergence of the generators implies a rate of convergence of the corresponding semigroups. This in turn will imply a rate of convergence in the vague topology of the law induced by $T_n$ to the one induced by $T$. Our result therefore establishes a sort of \emph{weak Berry-Esseen estimate}. The primary applications we have in mind are boundary problems, where we show that, in many cases, it is crucial that the function $f_n$ does not necessarily coincide with $\pi_n f$ but can be chosen to address issues arising from the boundary conditions. As an interesting application of this general framework, we establish a functional central limit theorem for a wide class of random walks on the positive integers, which converge to the most general Brownian motion on the positive half-line.

	The most general Brownian motion on the positive half-line was studied by Feller, and a comprehensive overview on it can be found in \citet[Theorem~6.2, p.\ 157]{Knight}. To put it simply, it is defined as a class of Feller processes on the positive half line such that its excursion to zero are the same as those of a standard Brownian motion, which can be shown to be a mixture of the reflected Brownian motion,  absorbed Brownian motion, and killed Brownian motion. Its generator is given by one half of the continuous Laplacian acting on the domain of $\mc C^2$-functions decaying at infinity and satisfying $c_1 f(0)-c_2 f'(0)+\frac{c_3}{2}f''(0)=0$ with $c_1+c_2+c_3=1$,  $c_i\geq 0$. Given three non-negative parameters $c_1, c_2, c_3$ that sum to one, we also occasionally denote by $B(c_1, c_2,c_3)$ the corresponding general Brownian motion.

	The discrete class of models considered here involves continuous-time random walks on $(\frac{1}{N}\mathbb{N})\cup\{\Delta\}$, where $\mathbb{N}={0,1,\ldots}$ and the state $\Delta$ is usually referred to as the \textit{cemetery}. The walk follows the usual symmetric walk on ${1,2,\ldots}$ with jump rates to nearest neighbors everywhere equal to $1/2$. However, at state $0$, we introduce the following: the rate to jump to state $1$ is $A/n^\alpha$, and the rate to jump to the cemetery is $B/n^\beta$, where $\alpha, A, \beta, B\geq 0$, and $n$ is the scaling parameter. Additionally, if the walk reaches the cemetery, it remains there indefinitely.
	
	The chosen values of parameters $\alpha, A, \beta, B$ then determine the limiting Brownian-type process of the random walk.
	We show here that for any choice of $c_1, c_2, c_3\geq 0$ such that $c_1+c_2+c_3=1$ and $c_1\neq 1$ there are classes of choices of $\alpha, A, \beta, B$ such that the above random walk converges to $B(c_1, c_2, c_3)$. We additionally show that by making a small shift to the right of the random walk,  we have convergence to the killed Brownian motion, which corresponds to  $c_1=1$.
	
	Each type of Brownian motion has different properties and behaviors, which makes them useful in different applications. For example, reflected Brownian motion can model a financial asset that cannot have negative values, while absorbed Brownian motion can model the extinction of a biological population. Elastic Brownian motion can model the behavior of a particle that is attracted to $0$ but has long-range repulsion from some boundary point and finally, the sticky Brownian motion can model a particle that sticks to a point, the killed Brownian motion corresponds to a continuous walk that jumps directly to the cemetery when it gets arbitrarily close to the origin. All these processes are particular cases of the above-mentioned general Brownian motion on the half-line, which illustrates the importance of the functional central limit theorem we prove in this paper. Moreover, the speed of convergence here obtained can come in handy when simulating those processes. Indeed, the random walks considered can be easily simulated as we will show in Section~\ref{s4}.
	
	The outline of the paper is the following: in Section~\ref{s2} we state results and mention related literature. Proofs are given in  Section~\ref{s3}. Finally, Section~\ref{s4} is dedicated to computational simulations of the simple random walk with boundary conditions which, by the presented results, are good approximations of the general BM on the half-line.

	\section{Statement of results}\label{s2}
	
	We use the convention that the set of natural numbers starts at zero, i.e., $\bb{N}=\{0,1,\ldots\}$.
	Throughout this paper, for given functions $f$ and $g$, the notation $f \lesssim g$ indicates the existence of a constant $c > 0$ such that $f(x) \leq cg(x)$ for all $x$ in the domain of the functions $f$ and $g$. We moreover use the notation $\bb R_{\geq 0}=[0,\infty)$.
	\subsection{The setup and main result}
	
	For each $n \in \bb{N}$, consider strongly continuous contraction semigroups $\{T_n(t) : t \geq 0\}$ and $\{T(t) : t \geq 0\}$ acting on Banach spaces $\msf{B}_n$ and $\msf{B}$, with generators $\msf L_n$ and $\msf L$, respectively. We denote the domains of $\msf{L}$ and $\msf{L}_n$ by $\mf{D} (\msf{L})$ and $\mf{D}(\msf{L}_n)$. To keep the notation simple, all norms are denoted by $\norm{\cdot}$.

	Let $\pi_n : \msf{B} \to \msf{B}_n$ be bounded linear operators indexed by $n \in \bb{N}$, which we call \textit{natural projections} for reasons that will become clear in the sequel.
	Let  $\Xi_n : \mf{D}(\msf{L}^2)\subset \msf{B} \to \msf{B}_n$ be a family of linear operators (not necessarily bounded), which we call \textit{the correction operators}. Denote by $\norm{\cdot}_\text{op}$ the operator norm. We further denote
	\begin{equation}\label{eq:A}
		\Phi_n \;=\; \pi_n + \Xi_n\,.
	\end{equation}
	The following hypotheses will be crucial in this work. Here, by $f \in \mf{D}(\msf{L}^2)$ we mean that $f \in \mf{D}(\msf{L})$ and $\msf{L}f \in \mf{D}(\msf{L})$.
	\begin{enumerate}[label=(H\arabic*)]
		\item \label{(H1)} If $f\in \mf{D}(\msf{L})$, then $\Phi_nf\in\mf D(\msf{L_n})$.
		\item \label{(H2)} There exist sequences $s_1(n)\searrow0$, $s_2(n)\searrow0$ and $s_3(n)\searrow 0$ such that, for any $f\in\mf{D}(\msf{L}^2)$,
		\begin{equation*}
			\norm{\pi_n\msf{L}f-\msf{L}_n\Phi_nf}\;\leq\;s_1(n)\norm{f}+s_2(n)\norm{\msf{L}f}+s_3(n)\norm{\msf{L}^{2}f}\,.
		\end{equation*}
		\item \label{(H3)} There exist sequences $r_1(n)\searrow0$ and $r_2(n)\searrow0$ such that, for each $f\in\mf{D}(\msf{L})$,
		\begin{equation*}
			\norm{\Xi_nf}\;\leq\;r_1(n)\norm{f} +r_2(n)\norm{\msf{L}f}\,.
		\end{equation*}
		\item \label{(H4)} For each $f\in \mf{D}(\msf{L}^2)$, the function $s \in [0,t] \mapsto \Phi_n T(t-s)f\in \msf{B}_n$ is differentiable and its derivative is given by $-\Phi_n \msf{L} T(t-s) f$.
	\end{enumerate}
	\begin{remark}\label{rmk21}
		Since $\pi_n$ is a linear bounded operator, it is easy to check that   $s \in [0,t] \mapsto \pi_n T(t-s)f\in \msf{B}_n$ is differentiable with derivative  $-\pi_n \msf{L}T(t-s) f$. The key point in \ref{(H4)} is to assure that    $s \in [0,t] \mapsto \Xi_n T(t-s)f\in \msf{B}_n$ is differentiable with derivative  $-\Xi_n \msf{L} T(t-s) f$.
	\end{remark}

	As we see next, these conditions  yield a bound on the distance between the semigroups generated by $\msf L_n$ and $\msf L$.
	This is a quantitative version of \citet[Theorem 6.1, p. 28]{EK}, attributed  to Sova and Kurtz,  which  is sometimes referred to  in the literature as \textit{Trotter-Kato Theorem}.
	\begin{theorem}\label{thm21}
		Under hypotheses  \ref{(H1)} -- \ref{(H4)}, for any $f \in \mf{D}(\msf{L}^2)$ and for each $t$ in a compact interval $[0,b]$ we have that
		\begin{align*}
			\norm{{T}_n(t)\pi_nf - \pi_nT(t)f}\;\lesssim\; & \max\left\{s_1(n),r_1(n)\right\}\norm{f} +\max\big\{s_2(n),r_1(n),r_2(n)\big\}\norm{\msf{L}f}\\
			&+\max\big\{s_3(n),r_2(n)\big\}\norm{\msf{L}^2f}\,.
		\end{align*}
		Here, the proportionality constant depends only on $b$.
	\end{theorem}

	Now we restrict ourselves to the setup of probability measures.
	Let $(\msf S,d)$ be a separable  locally compact (but possibly not complete) metric space and consider the one-point compactification $S_\Delta=S\cup \{\Delta\}$, where $\Delta$ is interpreted as cemetery point. We refer the reader also to~\cite{EK} for more on the one-point compactification.
	\begin{definition}\label{def:C_0}
		Let $\mc{C}_0(\msf S_\Delta)$ be the space of continuous functions $f:\msf S_\Delta\to \bb R$ such that $f(\Delta) =0$.
	\end{definition}
	\begin{remark}\label{rem:onepoint}
		The choice $f(\Delta)=0$ in the definition above is a convention. The cemetery $\Delta$ can be understood as the point that lies out of any compact set on $S$. 
		In this work the two relevant cases are $S=[0,\infty)$ and $S=(0,\infty)$. In the former case the space $\mc{C}_0(\msf S_\Delta)$ consists of all continuous functions converging to zero at infinity, and in the latter of all continuous functions converging to zero at zero and at infinity.
	\end{remark}
	We state below our next hypotheses:
	\begin{enumerate}[label=(G\arabic*)]
		\item \label{(G1)} The Banach space $\msf B$ is given by $\mc{C}_0(\msf S_\Delta)$ equipped with the uniform topology and the natural projection $\pi_n$ is the restriction to a subset $\msf S_n$ of $\msf S$, i.e., $\pi_n f= f\vert_{\msf S_n}$.
		\item \label{(G2)} There is a sequence of functions  $\{f_{k}\}_{k\geq 0}$ in $\mc C_0(\msf S_\Delta)$ such that $\spann(\{f_k\})$ is dense. Moreover, for each $k \in \bb{N}$, there exists a sequence of functions $\{f_{k,j}:j \geq 1\} \subset \mf{D}(\msf{L}^2)$ such that $f_{k,j} \to f_k$ when $j \to \infty$, in the uniform topology.
		\item \label{(G3)} 
		It holds that
		\begin{equation}\label{eq:fkj}
			\sum_{j\geq 1,k \geq 0}\frac{\norm{f_{k,j}}}{2^{k+j}} < \infty\;,\qquad \sum_{j\geq 1,k \geq 0}\frac{\norm{\msf{L}f_{k,j}}}{2^{k+j}} < \infty\,,\qquad\text{and}\qquad
			\sum_{j\geq 1,k \geq 0}\frac{\norm{\msf{L}^2f_{k,j}}}{2^{k+j}} < \infty\,.
		\end{equation}

		\item \label{(G4)} The semigroup $\msf{T}$ associated to the generator $\msf{L}$ is locally Lipschitz in the sense that, for each $t>0$ and each compact set $K\subseteq \msf{S}$, there exists a constant $M = M(t,K)>0$ such that for all $f\in \msf B=\mc{C}_0(\msf S)$
		\begin{equation*}
			\vert  \msf{T}(t)f (u) - \msf{T}(t)f (v) \vert \; \leq \; M\cdot \Vert f\Vert \cdot  d(u,v)\,,\quad \forall\, u,v\in K\,.
		\end{equation*}

	\end{enumerate}

	We now  introduce a suitable distance
	on the space of sub-probability measures that implies vague convergence.
	Let $\mc{B}$ be the Borel $\sigma$-algebra on $(\msf{S},d)$ and let $\mc M_{\leq 1}(\msf{S})$ be the set of sub-probability measures over the measurable space $(\msf{S}, \mc{B})$.
	\begin{definition}\label{def:metric}
		Let $\{f_{k}\}_{k \geq 0}\subset\mc{C}_0(\msf S_\Delta)$ and $\{f_{k,j}:j \geq 0\} \subset \mf{D}(\msf{L}^2)$ be as in \ref{(G2)}. We define
		\begin{equation*}
			{\bf d}(\mu,\nu)\;:=\; \sum_{j\geq 1,k=0}^\infty \frac{1}{2^{k+j}}\left(\bigg\vert \int f_{k,j}\dd \mu - \int f_{k,j}\dd \nu\bigg\vert \wedge 1\right)
		\end{equation*}
		for any $\mu, \nu\in \mc M_{\leq 1}(S)$.
	\end{definition}
	
	At a first glance, it may seem odd to define the distance ${\bf d}$ in terms of a double indexed sequence $\{f_{k,j}\}$ instead of a single indexed sequence, as it is usually done. We use the above notation since it will be more natural for our applications.
	We also highlight the fact that each $f_{k,j}$ is in the domain of $\msf L^2$, which will be crucial for the arguments to come.
	The next proposition is a straightforward consequence of the definitions of vague and weak convergence, and for this reason we omit its proof.
	\begin{proposition}
		The function ${\bf d}: \mc{M}_{\leq 1}(\msf S) \times \mc{M}_{\leq 1}(\msf S)\to [0,+\infty)$ defined in Definition~\ref{def:metric} is a metric on the set
		$\mc{M}_{\leq 1}(\msf S)$ and convergence with respect to  the metric $\bf d$ is equivalent to vague convergence. Moreover, if $\lim_{n\to\infty}{\bf d}(\mu_n,\mu)=0$ and  $\mu$ is a probability measure, then $\mu_n$ converges weakly to $\mu$.
	\end{proposition}

	We are in a position to state the first main result of this work. Recall that $\msf S_n\subset \msf S$. In the following theorem we use the notion of Feller processes which can be found in~\cite[Chapter 6.2, Definition~6.7]{LeGall}.
	\begin{theorem}[Berry-Esseen estimate with respect to ${\bf d}$ and functional CLT]\label{thm2.3}\quad \\
		Assume hypotheses \ref{(H1)}--\ref{(H4)} and \ref{(G1)}--\ref{(G4)}. Let  $\big\{X(t):t\in [0,T]\big\}$ and $\big\{X_n(t):t\in [0,T]\big\}$ be the Feller  processes on $\msf S$  and $\msf S_n$ associated to the generators $\msf{L}$ and $\msf{L}_n$,  respectively,   assumed to start at the points $x_n \in \msf{S_n}$ and  $x\in \msf{S}$, respectively, where $d(x_n,x)\leq i(n)$ for some function $i(n)$ converging to zero as $n\to\infty$.
		Fix some time  $t>0$ and denote by $\mu$ and $\mu_n$ the probability distributions on $\msf S$ induced by $X(t)$ and $X_n(t)$, respectively, starting from the points $x_n\in \msf{S}_n$ and $x\in \msf{S}$.  Then
		\begin{equation}\label{berryessen}
			{\bf d}(\mu,\mu_n)\;\lesssim \;\max\big\{i(n), s_1(n),s_2(n), s_3(n),r_1(n),r_2(n)\big\}\,,
		\end{equation}
		and weak convergence  $X_n\Rightarrow X$ holds in the Skorohod space $\msf{D}_{\msf S}[0, \infty)$.
	\end{theorem}

	\begin{remark}\rm
		The main novelty in the statement above is the weak Berry-Esseen estimate \eqref{berryessen}. The weak convergence in the Skorohod space $\msf{D}_{\msf S}[0, \infty)$, denoted by $\Rightarrow$, is actually a corollary of Theorem~\ref{thm21} using \citet[Theorem 2.11, page 172]{EK}. We refer also to~\cite[Chapter 3]{EK} for the definitions of weak convergence and the Skorohod space $\msf{D}_{\msf S}[0, \infty)$.
	\end{remark}

	\subsection{Application to the  general Brownian motion on the half-line}
	
	We apply Theorem~\ref{thm21} in order to obtain a Donsker-type theorem for what is known in the literature as the \textit{general Brownian motion on the half line}. To that end, let us first define it. We denote by $B_0(t)$ the Brownian motion on $[0,\infty)$ absorbed upon reaching zero. We use $\Delta$ to represent the cemetery state, and we let $T_0$ be the hitting time of zero.
	\begin{definition}[\citealp{Knight}, p.\ 153]
		\label{def:GBM}
		A general Brownian motion  on the positive half-line is a diffusion process $W$ on the set $\bb{G} := \{\Delta\}\cup \bb R_{\geq 0}$ such that the absorbed process $\{W(t\wedge T_0):t\geq 0\}$ on $[0,\infty)$ has the same distribution as $B_0$ for any starting point $x\geq 0$.
	\end{definition}
	In the following definition we use the notation $\mc{C}_0^2(\bb{G})$ for the set of functions $f\in \mc{C}_0(\bb{G})$ that are twice differentiable for all $x> 0$ such that the second derivative on $(0,\infty)$ can be extended to an element of $\mc{C}_0(\bb{G})$.
	\begin{theorem}[Due to Feller, see \citealp{Knight}, Theorem 6.2, p.\ 157]\label{dgBM}
		Any general Brownian motion $ W$ on $[0, \infty)$ has generator $\msf L=\frac{1}{2}\frac{d^2}{dx^2}$ with corresponding domain
		\begin{equation}\label{domainBM}
			\mf{D}\big(\msf{L}\big)=\Big\{f\,\in\,\mc{C}_0^2(\bb{G})\,:f''\in \mc{C}_0(\bb{[0,\infty)})\, \text{ and } c_1f(0)-c_2f'(0)+\frac{c_3}{2}f''(0)=0\Big\}
		\end{equation}
		for some $c_i\geq 0$ such that $c_1 + c_2 + c_3 = 1$ and $c_1 \neq 1$.
	\end{theorem}
	We will now discuss the above process.
	The case $c_2 = 1$ corresponds to the \emph{reflected Brownian motion} (which has the distribution of the modulus of a standard BM), while $c_3 = 1$ yields the absorbed Brownian motion, which has the distribution of a standard BM stopped at zero.
	\begin{figure}[!b]
		\centering
		\begin{tikzpicture}[scale=1]
			
			\coordinate (A) at (-1,0);
			\coordinate (B) at (1,3.46);
			\coordinate (C) at (3,0);
			\coordinate (E) at (1,1);
			
			\draw (A) -- node[midway, above, rotate=60, align = center] {\small{elastic BM}\\ $c_3=0$} (B);
			\draw (B) -- node[midway, above, rotate=-60, align = center] {\small{sticky BM}\\ $c_1=0$} (C);
			\draw (C) -- node[midway, below, align = center] {\small{exponential holding BM}\\ $c_2=0$} (A);
			\draw[very thick] (A) -- (B) -- (C) -- cycle;
			\fill[gray!20] (A) -- (B) -- (C) -- cycle;
			\draw (E)  node[align = center]{\small{mixed BM} \\$c_1, c_2, c_3>0$};
			\draw (A) node[left, align = center]{\small{killed BM}\\ $c_1=1$};
			\draw (B)+(0,0.5) node[align=center]{\small{reflected BM}\\$c_2=1$};
			\draw (C) node[right, align = center]{\small{absorbed BM}\\ $c_3=1$};
			\filldraw[fill=green] (A) circle (3pt);
			\filldraw[fill=blue] (B) circle (3pt);
			\filldraw[fill=red] (C) circle (3pt);
		\end{tikzpicture}
		\caption{Description of the general Brownian motion on the half-line according to the chosen values of $c_1, c_2, c_3\geq 0$ on the simplex $c_1+c_2+c_3=1$. Note that the killed BM, which formally corresponds to $c_1=1$,   is not rigorously a case of the general BM, as we can see in Theorem~\ref{dgBM}.}
		\label{FigSimplex}
	\end{figure}
	As one can see above, the case $c_1=1$ is not part of Feller's Theorem~\ref{dgBM}, which imposes $c_1\neq 1$.   Indeed, for $c_1=1$ the domain \eqref{domainBM} is not a dense set in $\mc C_0(\bb G)$, where $\bb G = \{\Delta\}\cup \bb R_{\geq 0}$, so it cannot be a domain of a generator. But once we remove the origin, considering $\bb G = \{\Delta\}\cup (0,\infty)$ instead, it will define a Feller process because now the set of test functions is assumed to converge to zero at the origin (c.f.\ Definition~\ref{def:C_0} and refer to  \citealp[Chapter 2]{Chung} for details).
	The killed BM can be understood as the process  which jumps immediately to the cemetery $\Delta$ upon ``reaching the origin''. Actually, it never touches the origin, but gets arbitrarily close,  which explains why $c_1=1$ is not included in Theorem~\ref{dgBM}: since the killed BM does not touch zero, it cannot satisfy the condition of Definition~\ref{def:GBM}. Despite the absorbed BM and killed BM being different processes, they are similar in nature, and to improve the presentation we will freely consider the killed BM as being a particular case of the general BM on the positive half-line.

	The case $c_1=0$ corresponds to the \textit{sticky Brownian motion}, which is an interpolation between the absorbed BM and the reflected BM, spending a time of positive Lebesgue measure at zero (but not staying there for any non degenerated time interval). See \cite{Borodin} and \cite{Warren} on the subject, for example.

	The case $c_3 = 0$ corresponds to the \textit{elastic BM} (also called  \textit{partially reflected Brownian motion}) which it is a mixture of the reflected BM and killed BM. It  can be also constructed in terms of the local time at zero: we generate an exponential random variable $\tau$ (whose parameter is related to $c_1$ and $c_2$) and a path of the reflected BM. Once the local time at zero of the reflected BM reaches the value $\tau$, the process goes to the cemetery and stays there forever, see \cite{Lejay} and references therein.

	\begin{figure}[!htb]
		\centering
		\begin{tikzpicture}[scale=1]
			\draw (1,0)--(8,0);
			
			\draw (4.5,0.6) node[above]{$\displaystyle\frac{1}{2}$};
			\draw (5.5,0.6) node[above]{$\displaystyle\frac{1}{2}$};
			\draw (1.5,0.6) node[above]{$\displaystyle\frac{B}{n^{\beta}}$};
			\draw (0,0.6) node[above]{$\displaystyle\frac{A}{n^{\alpha}}$};
			
			\centerarc[thick,<-](1.5,0)(30:140:0.6);
			\centerarc[thick,->](0,-0.9)(55:125:1.5);
			\centerarc[thick,->](4.5,0)(40:150:0.6);
			\centerarc[thick,<-](5.5,0)(30:140:0.6);
			
			\draw (-1,-0.2) node[below]{$\Delta$};
			\draw (1,-0.2) node[below]{$\frac{0}{n}$};
			\draw (2,-0.2) node[below]{$\frac{1}{n}$};
			\draw (3,-0.2) node[below]{$\frac{2}{n}$};
			\draw (4,-0.2) node[below]{$\frac{3}{n}$};
			\draw (5,-0.2) node[below]{$\frac{4}{n}$};
			\draw (6,-0.2) node[below]{$\frac{5}{n}$};
			\draw (7,-0.2) node[below]{$\frac{6}{n}$};
			\draw (8,-0.2) node[below]{$\frac{7}{n}$};
			\draw (8.7,0) node[left]{...};
			
			\filldraw[ball color=black] (1,0) circle (.15);
			\filldraw[fill=white] (2,0) circle (.15);
			\filldraw[fill=white] (3,0) circle (.15);
			\filldraw[fill=white] (4,0) circle (.15);
			\filldraw[fill=white] (5,0) circle (.15);
			\filldraw[fill=white] (6,0) circle (.15);
			\filldraw[fill=white] (7,0) circle (.15);
			\filldraw[fill=white] (8,0) circle (.15);
			\filldraw[fill=white] (-1,0) circle (.15);
		\end{tikzpicture}
		\caption{Jump rates for the \textit{boundary random walk.}}
		\label{GF1}
	\end{figure}
	
	The case $c_2 = 0$ corresponds to the  \textit{exponential holding BM}. Its behavior is the following: once it visits zero, it  stays there for an exponentially distributed amount of time and then is killed (jumps to the cemetery and stays there forever). See \cite{Knight}  on the exponential holding BM. This case allows us to interpret the case $c_1=1$ as a kind of explosion: it is an extreme case of the exponential holding BM, where the parameter of the exponential clock associated to jumps from the origin to the cemetery is  infinite, leading to an instantaneous jump.
	
	The case $c_1,c_2,c_3>0$ is a mixture of these behaviours, and we will refer to any such process as the \textit{mixed BM}.  Finally, it is instructive to mention that there are actually only two behaviours at zero. Namely, how much the BM sticks at zero, and how much the BM is attracted to the cemetery. This is in agreement with the fact that there are three parameters $c_1$, $c_2$ and $c_3$, but only two degrees of freedom, since those parameters  are restricted to the simplex $c_1+c_2+c_3=1$. See Figure~\ref{FigSimplex} for an illustration of the general BM in terms of the choices of $c_1$, $c_2$ and $c_3$.

	Next  we define a random walk on the non-negative integers that, as we shall see, is a discrete version of the general BM on the positive half-line.
	
	The \textit{boundary random walk} is the Feller process depending on parameters $A, B, \alpha,\beta\geq 0$ and $n\in \bb N$  on the state space $\bb{G}_n := \big(\frac{1}{n}\bb{N}\big) \cup \{\Delta\}$, whose  generator $\mathsf{L}_n$  acts on functions $f : \mathbb{G}_n \to \bb{R}$ that decay to zero at infinity and are such that $f(\Delta)=0$ as follows:
	\begin{equation}\label{RWgm}
		\msf{L}_nf(\pfrac{x}{n})=n^2\times
		\begin{cases}
			\dfrac{1}{2}\Big[f\big(\pfrac{x+1}{n}\big)-f(\pfrac{x}{n})\Big]+\dfrac{1}{2}\Big[f\big(\pfrac{x-1}{n}\big)-f(\pfrac{x}{n})\Big]\,,& x=1,2, \ldots  \vspace{3pt}\\
			\dfrac{A}{n^\alpha}\Big[f\big(\Delta\big)-f(\pfrac{0}{n})\Big]+\dfrac{B}{n^\beta}\Big[f\big(\pfrac{1}{n}\big)-f\big(\pfrac{0}{n}\big)\Big]\,,&  x= 0\,,\vspace{3pt}\\
			0\,,&x =\Delta\,.
		\end{cases}
	\end{equation}
	Note that the jump  rates are bounded from above, have finite range and therefore $\msf L_n$ is a bounded operator generating a Feller process, see \cite[Theorem~9]{Reuter}.
	The jump rates are illustrated in Figure~\ref{GF1}.  The next statement is our second main result, and as we will see is an application of  Theorem~\ref{thm2.3}. The different cases stated below are illustrated in Figure~\ref{fig1}.
	\begin{theorem}[Functional CLT for the boundary random walk]\label{thm26}
		Fix $u,t>0$.	Let $\{X_n(t): t\geq 0\}$ be the boundary random walk of parameters $\alpha, \beta, A, B\geq 0$ starting from the point $\frac{\lfloor un \rfloor}{n}\in \bb G_n\subset \bb G$  and denote by $\mu_n = \mu_n(t)$ the distribution of $X_n$ at time $t>0$. Recall the metric ${\bf d}$ defined in Definition~\ref{def:metric} and denote by $\mu=\mu(t)$ the distribution at time $t>0$ of the limit process in each of following cases.   Then:
		\begin{enumerate}[leftmargin=*]
			\item If $\alpha=\beta +1$ and $\beta\in[0,1)$, then $\{X_n(t): t\geq 0\}$ converges weakly to $\{X^\text{EBM}(t):t\geq 0\}$ in the Skorohod topology of $\msf{D}_{\bb{G}}[0, \infty)$, where $X^\text{EBM}$ is  the elastic BM on $\bb G= \{\Delta\}\cup \bb R_{\geq 0}$ of parameters
			$$c_1 = \frac{B}{A+B},\quad c_2= \frac{A}{A+B} \quad \text{ and }\quad  c_3=0$$
			starting from the point $u$. Moreover,\medskip
			\begin{enumerate}
				\item if $\beta \in (0,1)$, then ${\bf d}(\mu_n,\mu) \lesssim \max\{n^{-\beta},n^{\beta-1}\}$\,,
				\item if $\beta = 0$, then ${\bf d}(\mu_n,\mu) \lesssim n^{-1}$\,.
			\end{enumerate}
			\medskip
			\item If $\alpha \in (2,+\infty]$  and $\beta =1$, then $\{X_n(t): t\geq 0\}$ converges weakly to $\{X^\text{SBM}(t): t \geq 0\}$ in the Skorohod topology of $\msf{D}_{\bb G}[0,\infty)$, where $X^\text{SBM}$ is  the sticky BM  on $\bb R_{\geq 0}$ of parameters $$c_1 = 0,\quad c_2 = \frac{B}{B + 1}, \quad \text{ and }\quad c_3 = \frac{1}{B + 1}$$
			starting from the point $u$. Moreover, in this case,  ${\bf d}(\mu_n,\mu) \lesssim \max\{n^{2 - \alpha},n^{-1}\}$.
			\item If $\alpha =2$ and $\beta\in(1,+\infty]$,  then  $\{X_n(t): t\geq 0\}$ converges weakly to $ \{X^\text{EHBM}(t):t\geq 0\}$ in the Skorohod topology of $\msf{D}_\bb{G}[0, \infty)$, where $X^\text{EHBM}$ is  the exponential holding BM  on $\bb G= \{\Delta\}\cup \bb R_{\geq 0}$ of parameters $$c_1 = \frac{A}{1+A},\quad c_2 = 0 \quad \text{ and }\quad c_3 = \frac{1}{ 1+A}$$
			starting from the point $u$. Moreover, for $\beta \in (2, \infty)$, ${\bf d}(\mu_n,\mu) \lesssim \max\{n^{2 - \beta}, n^{-1}\}$.
			\item If $\alpha>\beta +1$ and $\beta\in[0,1)$, then $\{X_n(t): t\geq 0\}$ converges weakly to $\{X^\text{RBM}(t): t\geq 0\}$ in the Skorohod topology of $\msf{D}_{\bb G}[0, \infty)$, where $X^\text{RBM}$ is  the reflected BM on $\bb{R}_{\geq 0}$, of parameters $$c_1 = 0,\quad c_2= 1 \quad \text{ and }\quad c_3=0$$
			starting from the point $u$. Moreover, 
			${\bf d}(\mu_n,\mu)\lesssim \max\{n^{\beta-1},n^{-\alpha+\beta+1 }\}\,.$
			\medskip
			\item If $\alpha\in(2,\infty]$ and $\beta \in (1,+\infty]$, then  $\{X_n(t): t\geq 0\}$ converges weakly to $\{X^\text{ABM}(t):t\geq 0\}$ in the Skorohod topology of $\msf{D}_{\bb G}[0, \infty)$, where $X^\text{ABM}$ is  the absorbed BM on $\bb{R}_{\geq 0}$, of parameters $$c_1 = 0,\quad c_2= 0\quad \text{ and } \quad c_3=1$$
			starting from the point $u$. Moreover, for $\alpha >2$ and $\beta > 2$, ${\bf d}(\mu_n,\mu) \;\lesssim\; \max\{n^{2 - \alpha}, n^{2 - \beta}, n^{-1}\}$.
			\item If $\alpha =2$ and $\beta =1$, then $\{X_n(t): t\geq 0\}$ converges weakly to  $\{X^\text{MBM}(t):t\geq 0\}$ in the Skorohod topology of $\msf{D}_\bb{G}[0,\infty)$, where $X^\text{MBM}$ is  the mixed BM on $\bb{G}$ of parameters $$c_1 = \frac{A}{1+A+B},\quad c_2 = \frac{B}{1+A+B}\quad \text{ and }\quad c_3 = \frac{1}{1+A+B}$$
			starting from the point $u$. Moreover, ${\bf d}(\mu_n,\mu)\lesssim n^{-1}$.
			
		\end{enumerate}
	\end{theorem}
	
	\begin{figure}[!htb]
		\centering
		\begin{tikzpicture}[smooth, scale = 1]
			\begin{scope}[scale = 1, xshift=0cm]
				
				\fill[fill=lightgray] (1,2)--(4,2)--(4,4)--(1,4)--cycle;
				\fill[pattern=custom north west lines,hatchspread=8pt,hatchthickness=0.5pt,hatchcolor=black] (0,1)--(1,2)--(1,4)--(0,4)--cycle;
				
				\draw[line width = 1 mm , white] (0,1)--(1,2);
				\draw[ultra thick, dash pattern=on 3.5pt off 1.5pt, gray] (0,1)--(1,2);

				\draw[ultra thick, white] (1,2)--(4,2);
				\draw[ultra thick, dash pattern=on 3.5pt off 1.5pt] (1,2)--(4,2);

				
				\draw[ultra thick] (1,2)--(1,4);
				
				\filldraw[ball color=white, draw=black, thick] (1,2) circle (2.5pt);
				
				\draw[->] (-0.5,0)--(4.5,0) node[below]{$\beta$};
				\draw[->] (0,-0.5)--(0,4.5) node[left]{$\alpha$};
				\draw (0.1,1)--(-0.1,1) node[left]{$1$};
				\draw (0.1,2)--(-0.1,2) node[left]{$2$};
				\draw (0.1,4)--(-0.1,4) node[left]{$\infty$};
				\draw (1,0.1)--(1,-0.1) node[below]{$1$};
				\draw (2,0.1)--(2,-0.1) node[below]{$2$};
				\draw (4,0.1)--(4,-0.1) node[below]{$\infty$};
				\draw (0,0) node[anchor= north east]{$0$};
			\end{scope}


			\begin{scope}[yshift=-0.5cm, xshift=-1cm]

				\begin{scope}[xshift = 0cm, yshift = 0cm]
					
					\begin{scope}[yshift=-0.5cm, xshift = 0cm]
						\filldraw[ball color=white, draw=black, thick] (0.5,-0.25) circle (2.5pt);
						\draw (1,-0.25) node[right]{Mixed BM};
					\end{scope}
					\begin{scope}[yshift=-1.25cm]
						\fill[fill=lightgray] (0,0) rectangle (1,-0.5);
						\draw (1,-0.25) node[right]{Absorbed BM};
					\end{scope}
					\begin{scope}[yshift=-2cm]
						\fill[pattern=custom north west lines,hatchspread=8pt,hatchthickness=0.5pt,hatchcolor=black] (0,0) rectangle (1,-0.5);
						\draw (1,-0.25) node[right]{Reflected BM};
					\end{scope}
				\end{scope}

				\begin{scope}[xshift=3.5cm]

					\begin{scope}[yshift=-0.5cm]
						\draw[line width = 1.5pt,color=black] (0,-0.25)--(1,-0.25);
						\draw (1,-0.25) node[right, color=black]{Stick BM};
					\end{scope}
					\begin{scope}[yshift=-1.25cm]
						\draw[ultra thick, dash pattern=on 3.5pt off 1.5pt](0,-0.25)--(1,-0.25);
						\draw  (1,-0.25) node[right]{Exp. Holding BM};
					\end{scope}
					\begin{scope}[yshift=-2cm]
						\draw[ultra thick, dash pattern=on 3.5pt off 1.5pt, color=gray] (0,-0.25)--(1,-0.25);
						\draw (1,-0.25) node[right]{Elastic BM};
					\end{scope}
				\end{scope}
				
			\end{scope}

		\end{tikzpicture}
		\caption{Possible limits for the boundary random walk according to the ranges of $\alpha,\beta\in[0,\infty]$. Note that it includes the cases where $\alpha =\infty$ or $\beta=\infty$, which correspond to $A=0$ and $B=0$ respectively. Speed of convergence is provided for all choices of $\alpha$ and $\beta$, except for the strip   $1< \beta < 2$.}
		\label{fig1}
	\end{figure}

	Since the natural state space of the killed BM is $\bb G = \{\Delta\}\cup (0,\infty)$, which does not include the origin, we need a different setup to have the convergence of the boundary random walk towards the killed BM.
	This is the content of the next result.
	Let $\tau_n : \bb G \to \bb G$ be the shift to the right of $1/n$ given by
	$\tau_n (\Delta) = \Delta$  and $\tau_n(u) = u+\pfrac{1}{n}$ for $u\in (0,\infty)$.
	\begin{theorem}[Convergence of the shifted boundary RW to the killed BM]\label{thm:killed}
		If $\alpha<1+\beta$ for $\beta\in [0,1]$, or $\beta>1$, then $\{\tau_n X_n(t): t\geq 0\}$ converges weakly to  $\{X^\text{KBM}(t):t\geq 0\}$ in the Skorohod topology of $\msf{D}_\bb{G}[0,\infty)$, where in this case $X^\text{KBM}$ is  the killed BM on $\bb{G} = \{\Delta\}\cup(0,+\infty)$, which is formally the general BM of parameters $c_1 = 1$, $c_2 = 0$ and $c_3 = 0$.
	\end{theorem}
	
	\begin{figure}[t]
		\centering
		\begin{tikzpicture}[scale=1]
			\draw (1,0)--(8,0);
			
			\draw (4.5,0.6) node[above]{$\displaystyle\frac{1}{2}$};
			\draw (5.5,0.6) node[above]{$\displaystyle\frac{1}{2}$};
			\draw (1.5,0.6) node[above]{$\displaystyle\frac{B}{n^{\beta}}$};
			\draw (0,0.6) node[above]{$\displaystyle\frac{A}{n^{\alpha}}$};
			
			\centerarc[thick,<-](1.5,0)(30:140:0.6);
			\centerarc[thick,->](0,-0.9)(55:125:1.5);
			\centerarc[thick,->](4.5,0)(40:150:0.6);
			\centerarc[thick,<-](5.5,0)(30:140:0.6);
			
			\draw (-1,-0.2) node[below]{$\Delta$};
			\draw (1,-0.2) node[below]{$\frac{1}{n}$};
			\draw (2,-0.2) node[below]{$\frac{2}{n}$};
			\draw (3,-0.2) node[below]{$\frac{3}{n}$};
			\draw (4,-0.2) node[below]{$\frac{4}{n}$};
			\draw (5,-0.2) node[below]{$\frac{5}{n}$};
			\draw (6,-0.2) node[below]{$\frac{6}{n}$};
			\draw (7,-0.2) node[below]{$\frac{7}{n}$};
			\draw (8,-0.2) node[below]{$\frac{8}{n}$};
			\draw (8.7,0) node[left]{...};
			
			\filldraw[ball color=black] (1,0) circle (.15);
			\filldraw[fill=white] (2,0) circle (.15);
			\filldraw[fill=white] (3,0) circle (.15);
			\filldraw[fill=white] (4,0) circle (.15);
			\filldraw[fill=white] (5,0) circle (.15);
			\filldraw[fill=white] (6,0) circle (.15);
			\filldraw[fill=white] (7,0) circle (.15);
			\filldraw[fill=white] (8,0) circle (.15);
			\filldraw[fill=white] (-1,0) circle (.15);    \end{tikzpicture}
		\caption{Jump rates for the \textit{shifted} boundary random walk.}
		\label{GF2}
	\end{figure}
	
	\section{Proofs}\label{s3}
	\subsection{Proofs of Theorem~\ref{thm21} and Theorem~\ref{thm2.3}}
	
	
	\begin{proof}[Proof of Theorem \ref{thm21}]
		Putting together hypotheses \ref{(H2)} and \ref{(H3)},   we can infer that
		\begin{align}\label{eq3.1}
			\norm{\Phi_n\msf{L}f - \msf{L}_n\Phi_nf}\,&\leq\,s_1(n)\norm{f}+2\max\{s_2(n),r_1(n)\}\norm{\msf{L}f}+2\max\{s_3(n),r_2(n)\}\norm{\msf{L}^2f}\,,
		\end{align}
		for any $f \in \mf{D}(\msf{L}^2)$.
		Fix $t \geq 0$ and define $g_{s,t} := \msf{T}(t - s)f$ for $0 \leq s \leq t$. Then,
		\begin{equation}\label{derivative}
			\partial_s g_{s,t}\;=\;\partial_s\msf{T}(t\,-\,s)f\;=\;-\,\msf{L}\msf{T(t\,-\,s)f\;=\;-\,\msf{L}g_{s,t}}.
		\end{equation}
		Note that $\Phi_nf = \Phi_ng_{t,t}$. From \eqref{derivative}  we obtain that
		\begin{align*}
			\msf{T}_n(t)\Phi_nf\,&=\,\msf{T}_n(0)\Phi_ng_{0,t}\,+\,\int_0^t\partial_s(\msf{T}_n(s)\Phi_ng_{s,t})\,\dd s\\
			&\overset{\ref{(H4)}}{=}\,\msf{T}_n(0)\Phi_ng_{0,t}\,+\,\int_0^t\big[\msf{L}_n\msf{T}_n(s)\Phi_ng_{s,t}\,-\,\msf{T}_n(s)\Phi_n\msf{L}g_{s,t}\big]\,\dd s\,.
		\end{align*}
		Using the equality above, using that semigroups are contractions, that $T_n(0)$ is the identity mapping and using that  $\msf{T}(t)\msf{L} = \msf{L}\msf{T}(t) $ and $\msf{T}_n(t)\msf{L}_n = \msf{L}_n\msf{T}_n(t) $ for all $n \in \bb{N}$ and for all $t \geq 0$,  we have that
		\begin{align}
			\norm{\msf{T}_n(t)\Phi_nf-\Phi_n\msf{T}(t)f}\;=\;&\Big\Vert\int_0^t\msf{T}_n(s)\big[\msf{L}_n\Phi_ng_{s,t}-\Phi_n\msf{L}g_{s,t}\big]\,\dd s\Big\Vert\notag\\
			&\leq\;\int_0^t\norm{\msf{T}_n(s)}_\text{OP}\cdot \norm{\msf{L}_n\Phi_ng_{s,t}-\Phi_n\msf{L}g_{s,t}}\,\dd s\notag \\
			&\overset{\eqref{eq3.1}}{\leq}\;\int_0^t\Big[s_1(n)\norm{g_{s,t}}+2\max\{s_2(n),\,r_1(n)\}\norm{\msf{L}g_{s,t}}\Big]\,\dd s\notag \\ &\qquad+\int_0^t2\max\{s_3(n),\,r_2(n)\}\norm{\msf{L}^2g_{s,t}}\,\dd s
			\notag\\
			&=\!\int_0^t\Big[s_1(n)\norm{\msf{T}(t-s)f}+2\max\{s_2(n),r_1(n)\}\norm{\msf{L}\msf{T}(t - s)f}\Big]\dd s\notag \\ &\qquad+2\int_0^t\max\{s_3(n),r_2(n)\}\norm{\msf{L}^2\msf{T}(t - s)f}\,\dd s
			\notag\\
			&\leq\;2t\Big(s_1(n)\norm{f}+\max\{s_2(n),\,r_1(n)\}\norm{\msf{L}f}\Big)\notag \\
			&\qquad\,+\,2t\max\{s_3(n),\,r_2(n)\}\norm{\msf{L}^2f}\,,\label{last}
		\end{align}
		for every $f \in \mf{D}(\msf{L}^2)$. Note that  hypothesis \ref{(H1)} assures that $\msf{L}_n\Phi_ng_{s,t}$ above is well-defined. Recall that $\Phi_n = \pi_n + \Xi_n$. By the triangle inequality,
		\begin{align*}
			\norm{\msf{T}_n(t)\pi_nf-\pi_n\msf{T}(t)f}\,&\leq\,\norm{\msf{T}_n(t)\Phi_nf-\Phi_n\msf{T}(t)f}\,+\,\norm{\msf{T}_n(t)\Xi_nf}\,+\,\norm{\Xi_n\msf{T}(t)f}.
		\end{align*}
		Invoking hypothesis \ref{(H3)} and applying inequality \eqref{last} we  conclude the proof.
	\end{proof}

	\begin{proof}[Proof of Theorem \ref{thm2.3}]
		Let $\{f_k\}$ satisfy hypothesis \ref{(G2)}.
		Denote
		\begin{align*}
			a(n)  = \max\{r_1(n),s_1(n)\}\,,\quad	b(n)  = \max\big\{s_2(n), r_1(n), r_2(n)\big\} \quad \text{and}\quad
			c(n)  = \max\big\{s_3(n), r_2(n)\big\}\,.
		\end{align*}
		Thus
		\begin{align*}
			& {\bf d}(\mu_n, \mu) =\; \sum_{{j\geq 1},k \geq 0}\frac{1}{2^{k + j}}\left(\mud{\int f_{k,j}\,\dd\mu_n - \int f_{k,j}\,\dd\mu}\wedge 1\right)\\
			&\overset{\ref{(G4)}}{\lesssim} \;\max\Big\{i(n), \sum_{j\geq 1,k \geq 0}\frac{1}{2^{k + j}}\Big(\norm{\msf{T}_n(t)\pi_nf_{k,j}-\pi_n\msf{T}(t)f_{k,j}}\wedge 1\Big)\Big\}\\
			&\overset{\text{Thm. }\ref{thm21}}{\lesssim}\;\max\Big\{i(n),\sum_{{j\geq 1},k \geq 0}\frac{1}{2^{k + j}}\Big(a(n)\norm{f_{k,j}} + b(n)\norm{\msf{L}f_{k,j}} + c(n)\norm{\msf{L}^2f_{k,j}}\Big)\wedge 1\Big\}\\
			&\overset{\ref{(G3)}}{\lesssim} 		 \max\big\{i(n), a(n), b(n), c(n)\big\}\,.
		\end{align*}
		Finally, since $X_n(0)\Rightarrow X(0)$ and ${\bf d}(\mu_n, \mu)\to 0$ as $n\to\infty$,   by \citet[Theorem~2.11, page~172]{EK} we conclude that $X_n \Rightarrow X$ under the Skorohod topology.
	\end{proof}

	\subsection{Proof of Theorem~\ref{thm26}}
	\label{sec:3.2}
	
	In order to prove Theorem~\ref{thm26}, we will apply Theorem~\ref{thm2.3}, so we must verify  hypotheses \ref{(H1)}--\ref{(H4)} and \ref{(G1)}--\ref{(G4)}.
	We start with \ref{(G1)}--\ref{(G3)}, whose arguments are general, afterwards we will deal with \ref{(H1)}--\ref{(H4)} and \ref{(G4)}, which require specific arguments for each choice of the parameters. 
	Hypothesis \ref{(G1)} is immediate from the setup.

	Let us check hypotheses \ref{(G2)} and  \ref{(G3)}. We first construct the family $\{f_k\}_{k\geq 0}$. Recall that the domain of $\msf{L}$ is given by \eqref{domainBM}. To construct the desired families of functions we first define $\mc{A}_{\geq 0} := \big\{p(x)e^{-x^2}:  \text{where }p:\bb R_{\geq 0}\to \bb R  \text{ is a polynomial }\big\}$. 
	It is a well-known result  that the linear vector space generated by $\mc{A}_{\geq 0}$ is dense in $\mc {C}_0(\bb R_{\geq 0})$. However, we were not able to find a reference for that fact. For sake of completeness we include a proof here. Let
	\begin{equation*}
		\mc{A} := \Big\{p(x)e^{-x^2}:  p:\bb R\to \bb R  \text{ is a polynomial }\Big\}\,.
	\end{equation*}
	\begin{lemma}\label{lemma1} The set
		$\mc{A}$ is dense in $\mc{C}_0(\bb R)$.
	\end{lemma}
	\begin{proof}
		Suppose by contradiction that $\mc{A}$ is not dense. Then, by the Hahn-Banach Theorem, there exists a non-zero functional $\Lambda : \mc{C}_0(\bb R) \to \bb{R}$ such that $\Lambda\vert_{\mc{A}} \equiv 0$. By the Riesz-Markov Theorem, there exists a measure $\mu$ such that
		$\Lambda(f):=\int f \,\dd\mu$ 	for all $f\in\mc{C}_0(\bb R)$.
		By the Jordan decomposition theorem, there exist $\mu^+,\mu^-$ positive real-valued measures such that $\mu = \mu^+ - \mu^-$ where at least one of the two measures is finite.
		For any $g \in \mc{A}$, it follows that $g(x) = p(x)e^{-x^2}$ for some polynomial $p$, and $\Lambda(g) = 0$. Thus
		\begin{equation}\label{eqref1}
			\int p(x)e^{-x^2}\,\dd\mu^+(x)\;=\;\int p(x)e^{-x^2}\,\dd\mu^-(x)\,.
		\end{equation}
		Let $\nu^+(A) = \int_A e^{-x^2} \dd\mu^+$ and  $\nu^-(A) = \int_A e^{-x^2} \dd\mu^-$ be real-valued measures, and define $\rho^+, \rho^- : \bb{C} \to \bb{R}$ through $\rho^\pm(z) := \int e^{zx} \dd\nu^\pm(x)$. Observe that $\rho^\pm$ are well-defined since
		\begin{equation*}
			\mud{\int_{\bb{R}} e^{zx}\,\dd\nu^\pm(x)}\;=\;\mud{\int_{\bb{R}} e^{zx}e^{-x^2}\,\dd\mu^\pm(x)}\;\leq\;\int_{\bb{R}} \mud{e^{zx - x^2}}\,\dd\mu^\pm(x)<\infty\,.
		\end{equation*}
		To see that the last inequality is true, note that
		$\int e^{-x^2}\,\dd \mu(x)=0$.
		Since either $\mu^{+}$ or $\mu^{-}$ are finite we can conclude that
		$	\int e^{-x^2}\,\dd \mu^\pm(x)<\infty$,
		from which the claim follows.
		In view of  \eqref{eqref1}, we have that
		\begin{equation*}
			\frac{d^{(n)}}{dz^{(n)}}\rho^+(0)\; =\; \frac{d^{(n)}}{dz^{(n)}}\rho^-(0)
		\end{equation*}
		which lead to, by comparison of its power series, that $\rho^+(z) = \rho^-(z)$ for any $z \in \bb{C}$. Therefore, for all $s \in \bb{R}$, $\rho^+(is) = \rho^-(is)$ for all $s \in \bb{R}$, and thus $\nu^+ = \nu^-$, which guarantees that $\Lambda \equiv 0$, a  contradiction. Hence, $\mc{A}$ must be dense in $\mc{C}_0(\bb R)$.
	\end{proof}
	As a consequence of Lemma~\ref{lemma1}, we have
	\begin{corollary}\label{cor:3.2}
		The set 	$\mc{A}_{\geq 0}$ is dense in $\mc{C}_0(\bb R_{\geq 0})$.
	\end{corollary}

	Let  $\tilde f_k:\bb R_{\geq 0}\to \bb R$ be  defined by $\tilde f_k(x) := x^ke^{-x^2}$ for all $k\geq 0$, which  are illustrated in Figure~\ref{fig:fk}.
	By  Corollary~\ref{cor:3.2} we know that $\spann(\{\widetilde f_k\}_{k\geq 0})$ is dense in $\mc C_0(\bb R_{\geq 0})$ and an elementary calculation gives that $\Vert\widetilde f_0\Vert = 1$  and \begin{equation}
		\label{eq:maximum}
		\Vert\widetilde f_k\Vert =\big(\pfrac{k}{2}\big)^{\frac{k}{2}}e^{- \frac{k}{2}}\;\qquad \text{for all $k \geq 1 $.}
	\end{equation}
	Then, the family of functions defined by
	\begin{equation}\label{eq:f}
		f_k\;:=\;\frac{\widetilde f_k}{\Vert\widetilde f_k\Vert}\,,\quad \forall\, k\geq 0\,,
	\end{equation}
	still has the property that its span is dense in $\mc C_0(\bb R_{\geq 0})$.
	\begin{figure}[b]
		\centering
		\begin{tikzpicture}[scale=2.5,smooth];
			\draw[dashed] (1,0.37) -- (0,0.37) node[left]{$1/e$};
			\draw[black = solid,  very thick, color = black] plot   [domain=0:3.5](\x,{exp(-\x*\x )});
			\draw[black = solid,  very thick, color= magenta ] plot   [domain=0:3.5](\x,{\x*exp(-\x*\x )});
			\draw[black = solid,  very thick, color = blue] plot   [domain=0:3.5](\x,{\x*\x*exp(-\x*\x )});
			\draw[black = solid,  very thick, color = green] plot   [domain=0:3.5](\x,{\x*\x*\x*exp(-\x*\x )});
			\draw[black = solid,  very thick, color = yellow] plot   [domain=0:3.5](\x,{\x*\x*\x*\x*exp(-\x*\x )});
			\draw[dashed] (1,0.37) -- (1,0) node[below]{$1$};
			\draw (0,1) node[left]{$1$};
			\draw (0,0) node[anchor = north east]{$0$};
			\begin{scope}[xshift = 0.25cm]
				\draw[black = solid, very thick]  (3,1.25)-- (3.5,1.25) node[right]{$\widetilde{f}_0$};
				\draw[black = solid, very thick, color = magenta]  (3,1)-- (3.5,1) node[right, color = black]{$\widetilde{f}_1$};
				\draw[black = solid, very thick, color = blue]  (3,0.75)-- (3.5,0.75) node[right, color = black]{$\widetilde{f}_2$};
				\draw[black = solid, very thick,  color = green]  (3,0.5)-- (3.5,0.5) node[right, color = black]{$\widetilde{f}_3$};
				\draw[black = solid, very thick, color = yellow]  (3,0.25)-- (3.5,0.25) node[right, color = black]{$\widetilde{f}_4$};
			\end{scope}
			\draw[->] (0,-0.25)--(0,1.5);
			\draw[->] (-0.25,0)--(4,0) node[anchor=north]{$x$};
		\end{tikzpicture}
		\caption{Illustration of the functions $\widetilde{f}_{k}: \bb R_{\geq 0}\to \bb R$, $\widetilde{f}_k(x) := x^ke^{-x^2}$.}\label{fig:fk}
	\end{figure}
	
	Our goal now is to find  sequences $f_{k, j}\in \mf D(\msf{L}^2)$ fulfilling  hypotheses \ref{(G2)} and  \ref{(G3)}. Note that the  functions $f_k$ are smooth and,  for $k\geq 5$, the function itself and its first four derivatives at zero are zero. Recalling \eqref{domainBM}, this trivially implies that $f_k\in \mf D(\msf{L}^2)$ for $k\geq 5$. Therefore,  we define
	\begin{equation*}
		f_{k,j}\; := \; f_k \quad \text{ for } k\geq 5\,.
	\end{equation*}
	We first discuss how to verify \ref{(G3)} in this case. The first condition in \eqref{eq:fkj} is immediate since the $f_k$'s are normalized. The second condition in \eqref{eq:fkj} can be checked by noting that 
	$2\msf{L}\widetilde f_{k}(x) = k^{2} x^{k-2} e^{- x^{2}} - 4 k x^{k} e^{- x^{2}} - k x^{k-2} e^{- x^{2}}+ 4 x^{k+2} e^{- x^{2}} - 2 x^{k} e^{- x^{2}}
	$,
	bounding the maximum by the sum of the maxima and using~\eqref{eq:maximum}. The third condition in \eqref{eq:fkj} follows in the same way.
	We turn to the case $k\in\{1,2,3,4\}$. We define the  shift operator $\tau_{j}$ by
	\begin{equation*}
		\tau_{j}(f)(x) \;=\; f\Big(x - \pfrac{4}{j}\Big)\one_{[x \geq \frac{4}{j}]}(x)\,.
	\end{equation*}
	Define now an extension of $\tau_{j}(f_k)$ to the whole line through a reflection around the $y$-axis, that is,
	\begin{equation*}
		g_{k,j}(x) \;:=\; \begin{cases}
			(\tau_{j})(f_k)(x)\,, & \text{ if } x\geq 0\\
			(\tau_{j})(f_k)(-x)\,, & \text{ if } x< 0\\
		\end{cases} \quad \text{ for }k \in \{1,2,3,4\} \text{ and  }  {j\geq 1}\,,
	\end{equation*}
	which are continuous, but not smooth at the points $\pm 4/j$. To remedy this, consider the $\mc C^\infty$-approximation of identity  $\varphi_j : \bb{R} \to \bb{R}$ given by
	\begin{align*}
		\varphi_j(x) :=
		\begin{cases}
			\displaystyle\frac{1}{c_j}\exp{\left(- \frac{1}{1 - (jx)^2}\right)}\,, &\text{ if } \mud{jx} < 1\vspace{5pt}\\
			0\,, &\text{ otherwise}
		\end{cases}
	\end{align*}
	where $
	c_j := \int_{\bb{R}}\exp{\Big(- \frac{1}{1 - (jx)^2}\Big)}\one_{\{|jx|\leq 1\}}\dd x$
	is the normalizing constant. Note that $\varphi_j(x) = j\varphi_1(jx)$, which yields the following relation between the derivatives
	\begin{equation}\label{eq:scaling}
		\frac{d^{(i)}}{dx^{(i)}}\varphi_j(x) \;=\; j^{i+1}\frac{d^{(i)}}{dx^{(i)}}\varphi_1(jx)\,,\quad\forall\, {j\geq 1}\,.
	\end{equation}
	Define now $f_{k,j}(x) := (g_{k,j}\ast\varphi_j)(x)$ for $k \in \{1,2,3,4\}$ and  ${j\geq 1}$,
	which is smooth and a simple but tedious calculation shows that $f_{k,j} \to f_k$ uniformly as $j \to \infty$. Since $f_{k,j}$ is symmetric in a neighborhood around zero, we obtain that $f_{k,j} \in \mf{D}(\msf{L}^2)$ for any $k \in \{1,2,3,4\}$ and any $j\geq 0$.
	
	Denote $\norm{\frac{d^{(i)}}{dx^{(i)}}\varphi_1} = A_i$ for $i\in \bb N$. Observe that
	\begin{align*}
		\left| \frac{d^{(i)}}{dx^{(i)}}f_{k,j}(x)\right| &\;=\; \left|\frac{d^{(i)}}{dx^{(i)}}\int_{-\frac{1}{j}}^{\frac{1}{j}} g_{k,j}(x - y)\varphi_j(y)\,dy\right| \\
		&\;=\; \left|\int_{-\frac{1}{j}}^{\frac{1}{j}} g_{k,j}(x-y)\frac{d^{(i)}}{dx^{(i)}}\varphi_j(y)dy \right|\;\leq\;  2j^i\norm{g_{k,j}}A_i\,,
	\end{align*}
	where we made use of the scaling relation~\eqref{eq:scaling} to obtain the last inequality.
	Since $\norm{g_{k,j}} = \norm{f_k}=1$, we obtain that
	\begin{align*}
		\norm{\msf{L}f_{k,j}} \leq j^2A_2 \quad \text{ and }\quad
		\norm{\msf{L}^2f_{k,j}} \leq \frac12 j^4A_4\,,\quad\forall \,k\geq 1, {j\geq 1}\,.
	\end{align*}
	
	Now, it only remains to construct $f_{0,j}$ such that it verifies conditions \ref{(G2)} and  \ref{(G3)}. To that end, let $\msf{P}:\bb{R}_{\geq 0} \to \bb{R}$ be a polynomial such that
	both $e^{-x^2}+\msf{P}$ and $(e^{-x^2}+\msf{P})''$ satisfy the boundary condition in \eqref{domainBM} and such that additionally $\msf{P}(0)=0$.
	To continue, define the $\mc C^\infty$-bump function $b_1: \bb{R} \to \bb{R}_{\geq 0}$ by
	\begin{align*}
		b_1(t) \;=\; 1 -  \frac{\ell(t^2-1)}{\ell(t^2-1) + \ell(2-t^2)}
		\qquad \text{ where }
		\qquad
		\ell(t) \;=\;
		\begin{cases}
			e^{-\frac{1}{t}}, &\text{ if } t > 0\\
			0, &\text{ if } t \leq 0\,,
		\end{cases}
	\end{align*}
	and for $j\geq 1$ define $b_j(t)=b_1(jt)$.
	The function $0\leq b_j\leq 1$ is equal to one in an interval of size $1/j$ around the origin and zero outside the interval $[-\sqrt{2}/j,\sqrt{2}/j]$.
	Finally,  define
	\begin{equation*}
		f_{0,j}(x) \;:=\; e^{-x^2} + (b_j\msf{P})(x)\,.
	\end{equation*}
	The fact that $f_{0,j}$ vanishes for $x\geq \sqrt{2}/j$ and that $f_{0,j}(0)=1$ together with its continuity guarantee that $f_{0,j} \to e^{-x^2}$ as $j \to \infty$ in the uniform topology. Since the polynomial $e^{-x^2}+ \msf{P}$ satisfies the aforementioned boundary condition,  we also have $f_{0,j} \in \mf{D}(\msf{L}^2)$ for all positive integers $j$. Note that $b_j(x) = b_1(jx)$, yielding for all $i\geq 0$ and all $j\geq 1$
	\begin{equation*}
		\frac{d^{(i)}}{dx^{(i)}}b_j(x) \;=\; j^{i}\frac{d^{(i)}}{dx^{(i)}}b_1(jx)\,.
	\end{equation*}
	Since $\msf{L}f = \frac{1}{2}f ''$, it is immediate that
	\begin{equation}\label{eq98}
		|\msf{L}e^{-x^2}| \;\lesssim\; (1+x^2)e^{-x^2}\quad \text{ and }\quad |\msf{L}^2e^{-x^2}| \;\lesssim\; (1+x^4)e^{-x^2}\,.
	\end{equation}
	Now, since $b_j \equiv 0$ outside the compact set $[-\frac{\sqrt{2}}{j}, \frac{\sqrt{2}}{j}]$ and $\msf{P}$ is a polynomial, for all $i\geq 1$ one has that $\norm{\frac{d^{(i)}}{dx^{(i)}} (b_j\msf{P})} \lesssim j^i$.
	Hence, using equation~\eqref{eq98}, the product rule and the bound above, we can obtain upper bounds for the generator $\msf{L}$ and $\msf{L}^2$ norm
	$
	\norm{\msf{L}f_{0,j}} \lesssim \norm{\tilde{f}_0} + \norm{\tilde{f}_2} + j^2$
	and
	$\norm{\msf{L}^2f_{0, j}} \lesssim \norm{\tilde{f}_0} + \norm{\tilde{f}_4} + j^4
	$.
	Hence, we ensured that conditions \ref{(G2)} and  \ref{(G3)} are met.
	
	It is missing to verify  \ref{(G4)}, which will rely  on some knowledge about the semigroup of the limiting process.

	The limiting processes mentioned in Theorem~\ref{thm26} are the reflected, absorbed, mixed, sticky,  elastic, exponential holding and killed Brownian motion. Except the sticky BM  and the exponential holding BM, all corresponding semigroup's formulas  can be found in \citet[Appendix~1]{handbook}. However, the semigroup for the sticky BM on $[0,\infty)$ can be deduced by symmetry from the formula for the sticky BM on $\bb R$ \citep[page~123]{handbook}. We could not find the semigroup for the exponential holding BM in the literature, so we provide below its  transition probability:
	\begin{equation*}
		\begin{cases}
			p_t(\Delta, \{\Delta\}) = 1 &\\
			p_t(0, \{\Delta\}) = 1- e^{-\lambda t} & \\
			p_t(0, \{0\}) = e^{-\lambda t} & \\
			p_t(x, A) = p_t^{\text{absorbed BM}}(x,A)& \forall\,A\not\ni 0,\Delta\quad \forall\,x>0\\
			p_t(x, \{\Delta\}) = \int_{0}^t f_{\tau_0}(s,x) \big(1- e^{-\lambda (t-s)}\big)ds & \forall\,x>0\\
			p_t(x, \{0\}) = \int_{0}^t f_{\tau_0}(s,x) e^{-\lambda (t-s)}ds & \forall\,x>0\\	
		\end{cases}
	\end{equation*}
	for all $t\geq 0$, where $p_t^{\text{absorbed BM}}$ is the transition probability function of the absorbed BM and 
	\begin{equation*}
		f_{\tau_0}(t,x) \;=\; \frac{x}{\sqrt{2\pi t^3}}e^{-\frac{x^2}{2t}}\one_{[0,\infty)}(t)
	\end{equation*}
	is the density of the hitting time of zero for the absorbed BM.
	In the supplementary material of this paper, see \cite{efjp_supplement}, we provide details on the claim that \ref{(G4)} holds for all semigroups dealt with in this paper. 
	
	\begin{remark}
		\rm Hypothesis \ref{(G4)} has been used just once, in the proof of Theorem~\ref{thm2.3}, and its importance relies on the fact that the random walk and  its limiting process may not have the same starting point. For instance, in the setup of Theorem~\ref{thm26}, the boundary random walk starts from $\lfloor un\rfloor/n$, whereas its Brownian counterpart starts from $u>0$. If we assume that the discrete process and the limiting process start from the same point $u\in\msf{S}$,  hypothesis \ref{(G4)} can be dropped from Theorem~\ref{thm2.3}. This is possible in the setup of Theorem~\ref{thm26}, for instance, if we assume that the scaling parameter is given by $n = n(k) = 2^k$  and the initial point $u$ is a positive integer.
	\end{remark}

	We therefore must  verify the remaining assumptions.
	Those involve the correction operator $\Xi_n$, which is very model dependent. Before we study each model separately, we start with some generalities. We will always assume that $(\Xi_n f)(\Delta)=0$, and moreover all functions considered here satisfy $f(\Delta)=0$ and they have continuous derivatives up to order four. Thus,
	the generator of the random walk with boundary conditions  applied to the function  $\Phi_nf$ at zero is given by
	\begin{equation}\label{eq:3.8}
		\msf{L}_n\Phi_nf\pn{0}=-\frac{A}{n^{\alpha - 2}}f\pn{0}+\frac{B}{n^{\beta - 2}}\left[f\pn{1}-f\pn{0}\right]-\frac{A}{n^{\alpha - 2}}\Xi_nf(0)+\frac{B}{n^{\beta - 2}}\big[\Xi_nf\pn{1} - \Xi_nf(0)\big].
	\end{equation}
	Outside of zero, a Taylor expansion yields
	\begin{align}
		\msf{L}_n\Phi_nf\pn{x} \;=\; &  \frac12 f''\pn{x}+ \frac{1}{2\cdot 4!n^2}\left[f^{''''}(\theta) + f^{''''}(\eta)\right]  + \frac{n^2}{2}\Big[\Xi_nf\pn{x+1} + \Xi_nf\pn{x-1} - 2\Xi_nf\pn{x}\Big]\notag\\
		\;=\; &  \frac12\pi_n  f''\pn{x}+ \Vert  f''''\Vert \cdot  O(n^{-2})  + n^2\msf{L}_n \Xi_n f \pn{x}\,,\label{lastline}
	\end{align}
	for some $\theta \in (x/n, (x+1)/n)$ and $\xi \in ((x - 1)/n, x/n)$.

	\subsubsection{The elastic BM: the case  \texorpdfstring{$\beta \in [0,1), \alpha = \beta + 1$}{beta in [0,1), alpha=beta + 1}}
	Recall that in this case we set
	$c_1 = \frac{B}{A + B}$, $ c_2 = \frac{A}{A + B}$ and $c_3 = 0$.
	Denote the generator of the elastic Brownian motion by $\msf{L}_\text{EBM}$, whose  domain is given by
	$
	\mf{D}\left(\msf{L}_\text{EBM}\right)=\left\{f\,\in\,\mc{C}_0^2(\bb{G}):\frac{A}{A+B}f(0)-\frac{B}{A+B}f'(0)\;=\;0\right\}$.
	Let $f \in \mf{D}(\msf{L}^2_\text{EBM})$, which  implies the boundary conditions $Af(0) = Bf'(0)$ and $Af''(0) = Bf'''(0)$.
	Using this together with $\alpha = 1+\beta$ in Equation~\eqref{eq:3.8} and  a Taylor expansion, we get that
	\begin{align}
		\msf{L}_n\Phi_nf\pn{0}\;=\;& -\frac{A}{n^{\alpha - 2}}f\pn{0}+\frac{B}{n^{\beta - 2}}\left[\frac{f'(0)}{n}+\frac{f''(0)}{2!n^2}+ \frac{f'''(0)}{3!n^{3}} + \frac{f''''(\eta)}{4!n^{4}}\right]\notag\\
		&-\frac{A}{n^{\alpha - 2}}\Xi_nf(0)+\frac{B}{n^{\beta - 2}}\Big[\Xi_nf\pn{1} - \Xi_nf(0)\Big]\notag\\
		\;=\;&   \frac{B}{n^{\beta}}\left[\frac{f''(0)}{2!}+ \frac{A}{B}\frac{f''(0)}{ 3!n} + \frac{f''''(\eta)}{4!n^{2}}\right] \label{vanishes} \\
		& -\frac{A}{n^{\beta- 1}}\Xi_nf(0)+\frac{B}{n^{\beta - 2}}\Big[\Xi_nf\pn{1} - \Xi_nf(0)\Big]\label{sobra}
	\end{align}
	for some $0\leq \eta\leq 1/n$. Note that for $\beta>0$ the parcel \eqref{vanishes} is vanishing, and it is at this point where the correction operator enters the game.
	Assume first that $\beta \neq 0$; we will discuss the situation $\beta=0$ at the ending of this subsection.
	Define the correction operator in this case by
	\begin{equation}\label{eq:corr}
		\Xi_nf\pn{x} \; :=\; -\frac{\frac{1}{2}f''(0)}{An^{1 - \beta}\left(1 + \frac{1}{n}g\left(\frac{x}{n}\right)\right)}\,,
	\end{equation}
	where $g$ is some arbitrary  non-negative Lipschitz function of constant $K > 0$ satisfying $g(0) = 0 $ and $g(u) \to \infty$ as $u \to \infty$. Note that the  condition on the growth of $g$ is only necessary to assure that $\Xi_n f$ belongs to the domain of $\msf{L}_n$ which is given by
	$ \mf D(\msf{L}_n)  = \big\{f: \{\Delta\}\cup \pfrac{1}{n}\bb N \to \bb R  \text{ such that  }  \lim_{\frac{x}{n}\to \infty} f\pn{x} = 0 \text{ and } f(\Delta)=0\big\}$. 
	
	By Remark~\ref{rmk21}, in order to verify \ref{(H4)}, we only need to check that $s \mapsto \Xi_n T(t-s)f\in \msf{B}_n$ has derivative given by $-\Xi_n \msf{L}T_n(t-s)f$. This is true because by~\eqref{eq:corr} one has that for all $x$ 
	\begin{equation*}
		\Xi_n T_n(t-s)f\pn{x}\; =\; -\frac{\frac{1}{2}(T(t-s)f)''(0)}{An^{1 - \beta}\left(1 + \frac{1}{n}g\left(\frac{x}{n}\right)\right)}\,,
	\end{equation*}
	and $\frac{1}{2}(T(t-s)f)''(0) = \msf{L}T(t-s)f(0) = T(t-s)\msf{L}f(0)$. Here we are allowed to differentiate since $f\in \mf D(\msf{L}^2)$. The same argument holds for all subsequent cases, except for the killed BM, which will be treated separately. Nevertheless, in that specific case we do not apply the Theorem~\ref{thm21}, so \ref{(H4)} is not necessary.  From now on we do not mention \ref{(H4)} anymore, being understood that it holds whenever it is needed.
	
	We then have that
	\begin{align*}
		-\frac{A}{n^{\beta- 1}}\Xi_nf(0)+\frac{B}{n^{\beta - 2}}\Big[\Xi_nf\pn{1} - \Xi_nf(0)\Big]
		& = \frac12f''(0) + \frac{Bn}{A}\left[\frac{\frac{1}{n}g\left(\frac{1}{n}\right)}{1 + \frac{1}{n}g\left(\frac{1}{n}\right)}\right]\frac12f''(0)\\
		&=\; \msf{L}_{\text{EBM}}f(0) + \frac{BK}{A}\Vert \msf{L}_\text{EBM}f\Vert \cdot O(n^{-1})
	\end{align*}
	because $g(0) = 0$ and $g$ is $K$-Lipschitz. Plugging it into \eqref{vanishes} -- \eqref{sobra} yields
	\begin{equation}\label{EBM_0}
		\left|\left(\pi_n\msf{L}_{\text{EBM}}f-\msf{L}_n\Phi_nf \right)(0)\right|\,\lesssim\,\left(\frac{1}{n^\beta} +\frac{1}{n} +  \frac{1}{n^{1 + \beta}}\right)\norm{\msf{L}_{\text{EBM}}f} + \frac{1}{n^{2 + \beta}}\norm{\msf{L}^2_{\text{EBM}}f}\,.
	\end{equation}
	For $\frac{x}{n} \in \bb{G}_n \backslash \{0\}$, equation \eqref{lastline} indicates that we need to estimate
	\begin{align*}
		\msf{L}_n \Xi_nf\pn{x}
		&=\frac{\msf{L}_{\text{EBM}}f(0)n^2}{An^{1 - \beta}}\left[\frac{\frac{1}{n}\left(g\pn{x+1} - g\pn{x}\right)}{\left(1 + \frac{1}{n}g\left(\frac{x + 1}{n}\right)\right)\left(1 + \frac{1}{n}g\left(\frac{x}{n}\right)\right)}+\frac{\frac{1}{n}\left(g\left(\frac{x-1}{n}\right) - g\left(\frac{x}{n}\right)\right)}{\left(1 + \frac{1}{n}g\left(\frac{x - 1}{n}\right)\right)\left(1 + \frac{1}{n}g\left(\frac{x}{n}\right)\right)}\right]\\
		& = K\Vert \msf{L}_{\text{EBM}}f \Vert O(n^{\beta-1})\,,
	\end{align*}
	where we used again that  $g$ is $K$-Lipschitz. Plugging it into   \eqref{lastline}, we conclude that
	\begin{equation}\label{EBM_x}
		\left|\left(\pi_n\msf{L}_\text{EBM}f-\msf{L}_n\Phi_nf \right)\pn{x}\right|\;\lesssim\;\frac{1}{n^{1-\beta}}\norm{\msf{L}_{\text{EBM}}f} + \frac{1}{n^2}\norm{\msf{L}^2_{\text{EBM}}f}
	\end{equation}
	uniformly in $\frac{x}{n} \in \bb{G}_n \backslash\{0\}$.
	Putting together \eqref{EBM_0} and \eqref{EBM_x}, we infer that
	\begin{equation}\label{EBM1}
		\begin{split}
			\norm{\pi_n\msf{L}_{\text{EBM}}f -\msf{L}_n\Phi_nf } \;\lesssim\;& \max \Big\{\frac{1}{n^\beta},\frac{1}{n}, \frac{1}{n^{1-\beta}}\Big\} \norm{\msf{L}_{\text{EBM}}f}  +\max\Big\{\frac{1}{n^2}, \frac{1}{n^{2+\beta}}\Big\} \norm{\msf{L}^2_{\text{EBM}}f}\\
			=\; &   \max \Big\{\frac{1}{n^\beta}, \frac{1}{n^{1-\beta}}\Big\} \norm{\msf{L}_{\text{EBM}}f} + \frac{1}{n^2} \norm{\msf{L}^2_{\text{EBM}}f}\,.
		\end{split}
	\end{equation}
	In view of  \eqref{EBM1}, we have assured  hypothesis \ref{(H2)} and it is  only missing to check \ref{(H3)}. From \eqref{eq:corr}, we immediately get that
	$\norm{\Xi_n f} \lesssim \frac{1}{n^{1-\beta}} \norm{\msf{L}_{\text{EBM}}f}$,
	showing that \ref{(H3)} holds. Hence, Theorem \ref{thm2.3}  yields
	\begin{equation*}
		{\bf d} (\mu_n,\mu)\;\lesssim\;\max \left\{\frac{1}{n^\beta}, \frac{1}{n^{1-\beta}}, \frac{1}{n}, \frac{1}{n^2}\right\} \;=\; \max \left\{\frac{1}{n^\beta}, \frac{1}{n^{1-\beta}}\right\}
	\end{equation*}
	and that $\{X_{n}(t): t \geq 0\}$ weakly converges to $ \{X^\text{EBM}(t): t \geq 0\}$ under the Skorohod topology of $\msf{D}_\bb{G}[0,\infty)$. We thus can conclude this case.
	We come to the case $\beta=0$. In this case we define
	\begin{equation*}
		\Xi_nf\pn{x} \;:=\; -\frac{(1 - B)\frac{1}{2}f''(0)}{An\left(1 + \frac{1}{n}g\pn{x}\right)}\,.
	\end{equation*}
	Analogous arguments as above yield that 	${\bf d} (\mu_n,\mu)\lesssim \frac{1}{n}$ in this case. We omit the details.

	\subsubsection{The sticky BM: the case  \texorpdfstring{$\beta = 1, \alpha \in (2, \infty)$}{beta = 1, alpha in (2, infty)}}
	Recall that in this case we set
	$c_1 = 0$, $c_2 = \frac{B}{B + 1}$ and  $c_3 = \frac{1}{B + 1}$.
	We denote the generator of the sticky BM by $\msf{L}_\text{SBM}$, whose  domain is
	$\mf{D}\left(\msf{L}_\text{SBM}\right)=\left\{f\,\in\,\mc{C}_0^2(\bb{G}):-\frac{B}{B + 1}f'(0)+\frac{1}{2}\frac{1}{B + 1}f''(0)=0\right\}$.
	As we shall see, no correction will be necessary here, being enough to choose $\Xi_n \equiv 0$.
	Let $f \in \mf{D}(\msf{L}^2_\text{SBM})$, which  yields the boundary conditions $Bf'(0) = \frac{1}{2}f''(0)$ and $Bf'''(0) = \frac{1}{2}f''''(0)$.
	Keeping this in mind and also that  $\beta = 1$ and $\alpha \in (2, \infty)$, we obtain from Equation~\eqref{eq:3.8} that
	\begin{align*}
		\msf{L}_n\Phi_nf(0)& \;=\;-\frac{A}{n^{\alpha - 2}}f(0)+Bf'(0)\,+\,\frac{B}{2n}f''(0)+\frac{B}{3!n^2}f'''(0)+\frac{B}{4!n^3}f^{''''}(\eta)\\
		&\;=\;-\frac{A}{n^{\alpha - 2}}f(0)+\Big(\frac{1}{2}+\frac{B}{2n}\Big)f''(0)+\frac{1}{2\cdot 3!n^2}f''''(0)+\frac{B}{4!n^3}f^{''''}(\eta)
	\end{align*}
	for some $0\leq \eta\leq 1/n$.
	Thus,
	\begin{equation*}\label{SBM1.1}
		|\left(\pi_n\msf{L}_\text{SBM}f - \msf{L}_n\Phi_nf\right)(0)| \leq \frac{A}{n^{\alpha - 2}}\norm{f} + \frac{B}{n}\norm{\msf{L}_\text{SBM}f} + \left(\frac{B}{6n^2} + \frac{B}{6n^3}\right)\norm{\msf{L}_\text{SBM}^2f}\,.
	\end{equation*}
	Recalling~\eqref{lastline} and our choice of $\Xi_n$ yields
	\begin{equation}\label{SBM1}
		\norm{\pi_n\msf{L}_{\text{SBM}}f - \msf{L}_n\Phi_nf} \;\lesssim\; \frac{1}{n^{\alpha - 2}}\norm{f} + \frac{1}{n}\norm{\msf{L}_\text{SBM}f} + \frac{1}{n^2}\norm{\msf{L}_\text{SBM}^2f}\,.
	\end{equation}
	Thus, Theorem~\ref{thm2.3} shows that
	\begin{equation*}
		{\bf d}(\mu_n,\mu)\;\lesssim\;\max\left\{\frac{1}{n^{\alpha - 2}}, \frac{1}{n^2}, \frac{1}{n}\right\}\;=\; \max\left\{\frac{1}{n^{\alpha - 2}}, \frac{1}{n}\right\}\,,
	\end{equation*}
	and also that $\{X_{n}(t): t\geq 0\}$ converges weakly to $\{X^\text{SBM}(t): t \geq 0\}$ under the Skorohod topology of $\msf{D}_{\bb{R}_{\geq 0}}[0,\infty)$. Hence, we can conclude this case.

	\subsubsection{The exponential holding BM: the case  \texorpdfstring{$\alpha = 2$}{alpha = 2} and  \texorpdfstring{$\beta \in (1, \infty]$}{beta in (1, infty]}}
	Here we set
	$c_1 = \frac{A}{A + 1}$, $c_2 = 0$ and $c_3 = \frac{1}{A + 1}$.
	Denote by $\msf{L}_\text{EHBM}$ the generator of the exponential holding BM, whose domain is  $\mf{D}\left(\msf{L}_\text{EHBM}\right) = \left\{f \in \mc{C}_0^2(\bb{G}):\frac{A}{A + 1}f(0) + \frac{1}{2}\frac{1}{A + 1}f''(0) = 0\right\}$.
	Additionally to $\alpha = 2$ and $\beta \in (1, \infty)$, assume for now that $\beta>2$. The case $\beta \in(1,2]$ will be analyzed later.
	Consider the correction operator identically null, that is, $\Xi_n \equiv 0$. Let $f \in \mf{D}(\msf{L}^2_\text{EHBM})$, which  yields the boundary conditions $Af(0) = -\frac{1}{2}f''(0)$ and $Af''(0) = -\frac{1}{2}f''''(0)$. Then,
	\begin{equation}\label{eq:eq}
		\begin{split}
			\msf{L}_n\Phi_nf(0) & \;=\; - Af(0) + \frac{B}{n^{\beta - 2}}\Big[f\pn{1} - f(0)\Big] \;=\; \frac{1}{2}f''(0) + \frac{B}{n^{\beta - 2}}\Big[f\pn{1} - f(0)\Big]\,,
		\end{split}
	\end{equation}
	which implies $\left|\left(\pi_n\msf{L}_\text{EHBM}f - \msf{L}_n\Phi_nf\right)(0)\right|\; \leq\; \frac{2B}{n^{\beta - 2}}\norm{f}$.
	Consequently, taking~\eqref{lastline} into account,
	\begin{equation}\label{EHBM2}
		\norm{\pi_n\msf{L}_\text{EHBM}f - \msf{L}_n\Phi_nf}\;\lesssim\; \frac{1}{n^{\beta - 2}}\norm{f} + \frac{1}{n^2}\norm{\msf{L}_\text{EHBM}^2f}\,.
	\end{equation}
	Thus, Theorem~\ref{thm2.3} implies that
	\begin{equation*}
		{\bf d}(\mu, \mu_n) \;\lesssim\; \max\left\{\frac{1}{n^{\beta - 2}}, \frac{1}{n^2},\frac{1}{n}\right\}\;=\; \max\left\{\frac{1}{n^{\beta - 2}}, \frac{1}{n}\right\},
	\end{equation*}
	and also that $\{X_{n}(t): t\geq 0\}$ converges weakly to $\{X^\text{EHBM}(t):t \geq 0\}$ in the Skorohod topology of $\msf{D}_{\bb{G}}[0, \infty)$.\medskip

	As a consequence of~\eqref{eq:eq} and using that 
	$ n^{2-\beta}[f(1/n) - f(0)]\sim n^{1-\beta}f'(0)$	we get that\break
	$\norm{\pi_n \msf{L}_{\text{EHBM}} f-\msf{L}_n \pi_n f}\to 0$. By \citet[Proposition 1.8, page 53]{engel2000one} for any generator $\msf{L}$ associated to a strongly continuous semigroup, the set $\mf{D}(\msf{L}^2)$ is a core for $\msf{L}$. Thus, in particular $\mf{D}(\msf{L}_{\text{EHBM}}^2)$ is a core for $\msf{L}_{\text{EHBM}}$. Applying 
	\citet[Theorem~6.1, page~28 and Theorem~2.11, page~172]{EK} one can conclude the convergence towards the exponential holding BM. However, since the rate of convergence relies on the first derivative of $f$, we are not allowed to apply Theorem~\ref{thm2.3}, and no speed of convergence could be provided in this case.
	
	\subsubsection{The reflected BM: the case \texorpdfstring{$\beta \in [0,1), \alpha >  \beta+1$}{beta in [0,1),  alpha >  beta+1}}

	Denote by $\msf{L}_{\text{RBM}}$ the generator of the reflected Brownian motion, whose domain is
	$
	\mf{D}(\msf{L}_{\text{RBM}}) := \left\{f \in \mc{C}_0^2(\bb{G}): f'(0) = 0\right\}$.
	Let $f \in \mf{D}(\msf{L}_{\text{RBM}}^2)$. Thus $f'(0)  = f'''(0) =0$, which leads to
	\begin{equation}\label{eq:generator_ref}
		\begin{split}
			\msf{L}_n\Phi_nf(0) &= -An^{2 - \alpha}f(0) + \frac{B}{2n^\beta}f''(0) + \frac{B}{4!n^{2 + \beta}}f''''(\eta)\\
			&\qquad  - \frac{A}{n^{\alpha - 2}}\Xi_nf(0) +\frac{B}{n^{\beta - 2}}\left[\Xi_nf\pn{1} - \Xi_nf(0)\right] \,,
		\end{split}
	\end{equation}
	for some $0\leq \eta\leq 1/n$. 
	Let $\Xi_nf\pn{x} :=    \widehat{\Xi}_nf\pn{x} + \widetilde{\Xi}_nf\pn{x}$,
	where
	\begin{align*}
		\widehat{\Xi}_nf\pn{x} \;=\; \frac{Af(0)}{Bn^{\alpha - \beta-1}}h\pn{x} \quad \text{ and }\quad \widetilde{\Xi}_nf\pn{x} \;=\; \Big(-\frac{1}{n}+ \frac{1}{Bn^{1-\beta}}\Big)\frac{f''(0)}{2}h\pn{x}\,,
	\end{align*}
	and we assume that $h$ is a fixed smooth compactly supported function satisfying $h(0) = h''(0) = 0$, $h'(0)=1$.
	First, note that
	\begin{equation}\label{eq:Xi}
		\norm{\Xi_n f} \;\lesssim\; \frac{\norm{f}}{n^{\alpha-\beta-1}} + \frac{\norm{f''}}{n^{1-\beta}}
	\end{equation}
	which converges to zero since $\beta+1<\alpha<2$ and $\beta <1$, verifying hypothesis \ref{(H3)}. Our goal now is to check \ref{(H2)}.
	Since $h(0)=0$,
	\begin{equation*}
		\begin{aligned}
			&- An^{2-\alpha }\widehat{\Xi}_nf(0)	+Bn^{2-\beta}\left[\widehat{\Xi}_nf\pn{1} - \widehat{\Xi}_nf(0)\right]\\
			&\;=\; Bn^{2-\beta}\widehat{\Xi}_nf\pn{1} \;=\; A f(0)n^{3-\alpha}\left[h(0)+h'(0)\frac{1}{n} + \frac{h''(0)}{2!}\frac{1}{n^2}+  \frac{h'''(\theta)}{3!}\frac{1}{n^3}\right]\,,
		\end{aligned}
	\end{equation*}
	for some $\theta \in [0,1/n]$. Since $h(0) = h''(0) = 0$ and $h'(0)=1$, we conclude that
	\begin{equation*}
		\begin{aligned}
			- An^{2-\alpha }\widehat{\Xi}_nf(0)	+Bn^{2-\beta}\left[\widehat{\Xi}_nf\pn{1} - \widehat{\Xi}_nf(0)\right]\;=\;A f(0)n^{2-\alpha} + \norm{f}\cdot O(\pfrac{1}{n^\alpha})\,.
		\end{aligned}
	\end{equation*}
	On the other hand, using that $h(0)= h''(0)=0$ and $h'(0)=1$,
	\begin{equation*}
		\begin{aligned}
			&- An^{2-\alpha }\widetilde{\Xi}_nf(0)	+Bn^{2-\beta}\left[\widetilde{\Xi}_nf\pn{1} - \widetilde{\Xi}_nf(0)\right]
			\;=\; (-Bn^{1-\beta} + n)\frac{f''(0)}{2}h\pn{1}\\
			&\;=\;  (-Bn^{1-\beta} + n)\frac{f''(0)}{2} \left[h'(0)\frac{1}{n} + \frac{h''(0)}{2!}\frac{1}{n^2}+  \frac{h'''(\theta)}{3!}\frac{1}{n^3}\right]\;=\; 
			(-Bn^{1-\beta} + n)\frac{f''(0)}{2} \left[\frac{1}{n} +\frac{h'''(\theta)}{3!}\frac{1}{n^3}\right]
		\end{aligned}
	\end{equation*}
	for some $\theta \in [0,1/n]$.
	Therefore, recalling \eqref{eq:generator_ref},
	\begin{equation*}
		\begin{aligned}
			& |\pi_n\msf{L}_{\text{RBM}}f(0)-\msf{L}_n\Phi_nf(0)|
			\;\lesssim\; \frac{1}{n^\alpha}\|f\|+\frac{1}{n^2}\|\msf{L}_{\text{RBM}}f\|+\frac{1}{n^{2+\beta}}\|\msf{L}^2_{\text{RBM}}f\|\,.
		\end{aligned}
	\end{equation*}
	Let us deal with the convergence outside zero. Using a Taylor expansion on $h$  similarly as in~\eqref{lastline} it is easy to check that, for $\frac{x}{n} \in \bb{G}_n \backslash\{0\}$,
	\begin{align*}
		\big\vert \msf{L}_n\widehat{\Xi}_nf\pn{x}\big\vert \;\lesssim\; \frac{|f(0)|}{n^{\alpha - \beta-1}}\Big[\norm{h''} +  \frac{\norm{h''''}}{n^{2}}\Big] \;\lesssim\; \frac{\|f\|}{n^{\alpha-\beta-1}}\,.
	\end{align*}
	On the other hand, also for  $\frac{x}{n} \in \bb{G}_n \backslash\{0\}$, we have
	\begin{align*}
		\msf{L}_n\widetilde{\Xi}_nf\pn{x}=\norm{f''}O\big(\pfrac{1}{n^{1-\beta}}\big)\,.
	\end{align*}
	Putting together all those bounds with \eqref{eq:generator_ref} and \eqref{lastline}, we finally get
	\begin{equation}\label{RBM0.1}
		\begin{aligned}
			\norm{\pi_n\msf{L}_{\text{RBM}}f - \msf{L}_n\Phi_nf}\;\lesssim\;\max\left\{\frac{1}{n^{\alpha}}, \frac{1}{n^{\alpha-\beta-1}}\right\}\norm{f} + \frac{1}{n^{1-\beta}}\norm{\msf{L}_{\text{RBM}}f} + \max\left\{\frac{1}{n^{2 + \beta}},\frac{1}{n^2}\right\}\|\msf{L}^2_{\text{RBM}}f\|\,.
		\end{aligned}
	\end{equation}
	In view of \eqref{RBM0.1} and \eqref{eq:Xi}, we can apply Theorem~\ref{thm2.3}, hence giving us that
	\begin{equation*}
		\begin{aligned}
			{\bf d}(\mu_n,\mu) &  \;\lesssim\; \max\left\{\frac{1}{n},\frac{1}{n^\alpha},  \frac{1}{n^{\alpha - \beta-1}}, \frac{1}{n^{1 - \beta}}, \frac{1}{n^{2+\beta}}, \frac{1}{n^2} \right\}   \;=\; \max\left\{  \frac{1}{n^{\alpha-\beta -1}}, \frac{1}{n^{1-\beta}} \right\}
		\end{aligned}
	\end{equation*}
	and that  $\{X_{n}(t): t \geq 0\}$ weakly converges to $\{X^{\text{RBM}}(t): t\geq 0\}$ under the Skorohod topology on $\msf{D}_{\bb{R}_{\geq 0}}[0, \infty)$, finishing the proof.

	\subsubsection{The absorbed BM: the case \texorpdfstring{$\beta > 1, \alpha>2$}{beta>1, alpha >2}}
	Denote by $\msf{L}_\text{ABM}$ the generator of the absorbed Brownian motion, whose domain is given by
	$\mf{D}\left(\msf{L}_\text{ABM}\right) := \{f \in \mc{C}_0^2(\bb{G}):f''(0) = 0\}$.
	
	\subsubsection*{\textbf{Subcase \texorpdfstring{$\alpha > 2$}{alpha >2} and \texorpdfstring{$\beta > 2$}{beta>2}}}
	Consider a null correction $\Xi_n \equiv 0$. Then
	\begin{align*}
		\msf{L}_n\Phi_nf(0) &=-\frac{Af(0)}{n^{\alpha - 2}} + \frac{B}{n^{\beta - 2}}\left[f\pn{1} - f(0)\right]\,,
	\end{align*}
	thus combining this with~\eqref{lastline}
	\begin{equation}\label{ABM1}
		\norm{\pi_n\msf{L}_\text{ABM}f - \msf{L}_n\Phi_nf} \;\lesssim\; \max\left\{\frac{1}{n^{\alpha - 2}}, \frac{1}{n^{\beta - 2}}\right\}\norm{f} + \frac{1}{n^2}\|\msf{L}_\text{ABM}^2f\|\,.
	\end{equation}
	Denote by $\mu$ the probability measure of the absorbed BM at time $t>0$. In view of \eqref{ABM1}, we can apply Theorem~\ref{thm2.3}, which gives us
	\begin{equation*}
		{\bf d}(\mu_n, \mu) \;\lesssim\; \max\left\{\frac{1}{n}, \frac{1}{n^{\alpha - 2}}, \frac{1}{n^{\beta - 2}},\frac{1}{n^2}\right\}
		\;=\; \max\left\{\frac{1}{n}, \frac{1}{n^{\alpha - 2}}, \frac{1}{n^{\beta - 2}}\right\}
	\end{equation*}
	and $\{X_{n}(t): t \geq 0\}$ weakly converges to  $\{X^\text{ABM}(t):t \geq 0\}$ under the Skorohod topology of $\msf{D}_{\bb{R}_{\geq 0}}[0,\infty)$.
	
	\subsubsection*{\textbf{Subcase \texorpdfstring{$\alpha > 2$}{alpha >2} and \texorpdfstring{$\beta \in (1,2]$}{beta in (1,2]}}}
	Similarly to what happened in the exponential holding BM in the strip $1<\beta<2$, we can check that $\norm{\pi_n \msf{L}_{\text{ABM}} f-\msf{L}_n \pi_n f }\to 0$ for every $f \in \mf{D}(\msf{L}_{\text{ABM}}^2)$, so we can apply  \citet[Theorem~6.1, page 28 and Theorem~2.11, page~172]{EK}  to deduce  the convergence towards the Absorbed BM. However, since the rate of convergence relies on the first derivative of $f$, we cannot apply   Theorem~\ref{thm2.3}, hence no speed of convergence could be provided in this case.
	
	\subsubsection{The mixed BM: the case \texorpdfstring{$\alpha = 2$}{alpha =2} and \texorpdfstring{$\beta = 1$}{beta =1}}
	Consider
	$c_1 = \frac{A}{1 + A + B}$, $c_2 = \frac{B}{1 + A +B}$ and $c_3 = \frac{1}{1 + A + B}$.
	Let  $\msf{L}_\text{MBM}$ be the generator of the mixed Brownian motion, whose domain is given by
	\begin{equation*}
		\mf{D}\left(\msf{L}_\text{MBM}\right) \;:=\; \Big\{f \in \mc{C}_0^2(\bb{G}):\frac{Af(0)}{1 + A + B} - \frac{Bf'(0)}{1 + A +B} + \frac{\frac{1}{2}f''(0)}{(1 + A + B)} = 0\Big\}.
	\end{equation*}
	In this case, also no correction is needed, so define $\Xi_n\equiv 0$ to be the null operator.
	If $f \in \mf{D}(\msf{L}_{\text{MBM}}^2)$, then we get the boundary conditions $Bf'(0) = Af(0) + \frac{1}{2}f''(0)$ and $Bf'''(0) = Af''(0) + \frac{1}{2}f''''(0)$. Since $\alpha=2$ and $\beta=1$, we have that using~\eqref{eq:3.8}
	\begin{align*}
		\msf{L}_n\Phi_nf(0) &\;=\; \left(Bf'(0) - Af(0)\right) + \frac{B}{2!n}f''(0) + \frac{B}{3!n^2}f'''(0) + \frac{B}{4!n^3}f^{''''}(\eta)
	\end{align*}
	for some $0\leq \eta\leq 1/n$. Applying the boundary conditions, it implies that
	\begin{equation*}
		\Big|\big(\pi_n\msf{L}_\text{MBM}f - \msf{L}_n\Phi_nf\big)(0)\Big| \;\leq\; \left(\frac{B}{n} + \frac{A}{3n^2}\right)\norm{\msf{L}_\text{MBM}f} + \left(\frac{1}{3n^2} + \frac{B}{6n^3}\right)\|\msf{L}_\text{MBM}^2f\|\,.
	\end{equation*}
	By the above bound and \eqref{lastline}
	\begin{equation*}
		\norm{\pi_n\msf{L}_\text{MBM}f-\msf{L}_n\Phi_nf}\; \lesssim\; \frac{1}{n}\norm{\msf{L}_\text{MBM}f} + \frac{1}{n^2}\|\msf{L}_\text{MBM}^2f\|\,.
	\end{equation*}
	Denote by $\mu$ the probability measure of the mixed BM at time $t>0$.
	Thus, by  the previous inequality we can invoke Theorem~\ref{thm2.3}, concluding that
	\begin{equation*}
		{\bf d}(\mu_n, \mu) \;\lesssim\;  \max\left\{\frac{1}{n}, \frac{1}{n^2}\right\}\;=\; \frac{1}{n}
	\end{equation*}
	and that $\{X_{n}(t): t \geq 0\}$ weakly converges to $\{X^{\text{MBM}}(t): t\geq 0\}$ in the Skorohod topology on $\msf{D}_\bb{G}[0,\infty)$.

	\subsection{Proof of Theorem~\ref{thm:killed}}\label{subsec:killed} In this section we prove Theorem~\ref{thm:killed}. Before delving into the details we recall that the topology is different here, because $\msf{S} = (0,\infty)$. In this scenario, see Remark~\ref{rem:onepoint}, the functions in the space $\mc{C}_0(\msf S_\Delta)$ must converge to zero at zero, see  Definition~\ref{def:C_0}.
	Since the topology is different from the previous setup, we need to check again hypothesis \ref{(G2)}. Recall the definition of the functions $f_k$ in \eqref{eq:f}. We note that, for all $j,k\geq 1$, we have that $f_{k,j}(0)= f_{k,j}''(0)=0$, which implies $f_{k,j}\in \mf D(\msf{L}_{\text{KBM}})$ for all $j,k\geq 1$. In the proposition below we do not consider the function $f_0$ since $f_0(0)\neq 0$, so it does not belong to the
	space $\mc{C}_0\big((0,\infty)\big)$. 
	\begin{proposition}\label{prop:killeddom}
		Let $\mc{B} = \{f_k: k \geq 1\}$. Then $\spann(\mc{B})$ is  dense in $\mc{C}_0\big((0,\infty)\big)$.
	\end{proposition}
	\begin{proof}
		Immediate from Corollary~\ref{cor:3.2} and the fact that $f_k(0)=0$ for any  $k\geq 1$.
	\end{proof}
	Proposition \ref{prop:killeddom} together with its preceding comment assure \ref{(G2)}. To ease arguments, note  that $\tau_n X_n$ is the same as considering the boundary random walk $X_n$  on the state space in Figure~\ref{GF2}.
	Denote by $X^{\text{KBM}}$  the killed Brownian motion of coefficients $c_1=1$, $c_2=0$, $c_3=0$ and  let  $\msf{L}_\text{KBM}$ be its generator, whose domain  is given by
	$\mf{D}\left(\msf{L}_\text{KBM}\right) := \{f \in \mc{C}_0^2(\bb{G}):f(0) = 0\}$.
	Let $f\in \mf{D}(\msf{L}_\text{KBM}^2)$. By the definition of $\mc{C}_0(\bb{G})$ and since $\msf{L}_\text{KBM}f\in \mc{C}_0(\bb{G})$, we infer that  $f''(0)= 0$.
	By doing Taylor expansions around zero and using that $f(0)=f''(0)=0$, we get
	\begin{align}
		\msf{L}_nf\pn{1}
		&=\;  \frac{An^2}{n^\alpha}\Big[f(\Delta) - f\pn{1}\Big] +\frac{Bn^2}{n^\beta}\Big[ f\pn{2} - f\pn{1}\Big] \notag \\
		&=\; -An^{2-\alpha}\Big[f(0)+f'(0)\pfrac{1}{n} +\frac{f''(0)}{2!}\frac{1}{n^2}+\frac{f'''(\theta_1)}{3!} \frac{1}{n^3}\Big] \notag \\
		&+ Bn^{2-\beta}\Big[ f'(0)\pfrac{2}{n}+\pfrac{f''(0)}{2!}\big(\pfrac{2}{n}\big)^2 +\pfrac{f'''(\theta_2)}{3!}\big(\pfrac{2}{n}\big)^3- f'(0)\pfrac{1}{n}-\pfrac{f''(0)}{2!}\pfrac{1}{n^2} -\pfrac{f'''(\theta_3)}{3!}\pfrac{1}{n^3}\Big] \notag\\
		& = \; f'\pn{0}\Big[Bn^{1-\beta}-An^{1-\alpha}\Big] + O(\max\{\pfrac{1}{n^{1+\beta}},\pfrac{1}{n^{1+\alpha}}\})\,,\label{nonvanishing}
	\end{align}
	where $\theta_1, \theta_3\in [0, \frac{1}{n}]$ and $\theta_2\in [0,\frac2n]$.
	On the other hand, we have that assuming that our correction $\Xi_n$ is such that $\Xi_nf(\Delta)=0$,
	\begin{align}
		\msf{L}_n\Xi_nf\pn{1} &\; =\;  \frac{An^2}{n^\alpha}\Big[\Xi_nf(\Delta) - \Xi_nf\pn{1}\Big] +\frac{Bn^2}{n^\beta}\Big[ \Xi_nf\pn{2} - \Xi_nf\pn{1}\Big]\notag\\
		&\; =\;  -An^{2-\alpha}\Xi_nf\pn{1} +Bn^{2-\beta}\Big[ \Xi_nf\pn{2} - \Xi_nf\pn{1}\Big]\,.\label{nonvanishing_2}
	\end{align}
	Since $\msf{L}_{\text{KBM}}f(0) = \frac{1}{2}f''(0) = 0$, our goal is to find a correction $\Xi_n$ such that \eqref{nonvanishing_2} cancels the non-vanishing terms in \eqref{nonvanishing}.
	
	\subsubsection{Subcase \texorpdfstring{$\beta>1$}{beta>1}}
	In this case, the parcel $Bn^{1-\beta}$ in \eqref{nonvanishing} vanishes as $n\to\infty$, and we only need to deal with $An^{1-\alpha}$.
	Define
	\begin{align*}
		\Xi_nf\pn{x}\;:=\; -\frac{f'(0)\cdot n^{-1}}{1+\frac{g\pn{x-1}}{n}}
	\end{align*}
	for $\frac{x}{n} \geq \frac{1}{n}$,  where $g: \bb G \to \bb R$ is a non-negative compact supported smooth function such that $g(0)=0$. It is now straightforward to check that  $	\Vert  \pi_n \msf{L}_{\text{KBM}}f-\msf{L}_n\Phi_nf \Vert \to 0$.
	Thus,  applying  \citet[Theorem~6.1, page~28 and Theorem~2.11, page~172]{EK}, we conclude that $\{X_{n}(t): t \geq 0\}$ weakly converges to $\{X^{\text{KBM}}(t): t\geq 0\}$ under the Skorohod topology on $\msf{D}_\bb{G}[0,\infty)$. Note that $\norm{\Xi_n f} = \norm{f'}O(1/n)$, so we cannot apply our Theorem~\ref{thm2.3} and not speed of convergence is provided.
	
	\subsubsection{Subcase \texorpdfstring{$\beta\in [0,1], \alpha < 1+\beta$}{0=< beta =< 1, alpha < 1+beta}}
	In this scenario, the parcel $Bn^{1-\beta}$ does not vanish as $n\to\infty$, and we need to define an extra correction to deal with it.
	Define
	\begin{align*}
		\Xi_nf\pn{x}\;:=\; -\frac{f'(0)\cdot n^{-1}}{1+\frac{g\pn{x-1}}{n}} +
		\frac{\frac{B}{A}f'(0)\cdot n^{\alpha-\beta-1}}{1+\frac{g\pn{x-1}}{n}}
	\end{align*}
	for $\frac{x}{n} \geq \frac{1}{n}$ and $g$  is the same function as in the previous case. Note that the  condition $\alpha<1+\beta$ guarantees that the correction above goes to zero as $n\to\infty$.  It is now straightforward to check that  $	\Vert  \pi_n \msf{L}_{\text{KBM}}f-\msf{L}_n\Phi_nf \Vert \to 0$, leading us to conclude that $\{X_{n}(t): t \geq 0\}$ weakly converges to $\{X^{\text{KBM}}(t): t\geq 0\}$ under the Skorohod topology on $\msf{D}_\bb{G}[0,\infty)$.
	\begin{remark}\rm
		Note that the region corresponding to the parameters $(\beta,\alpha)$ for which the shifted boundary random walk converges to the killed BM covers
		the white region in the first quadrant of Figure~\ref{fig1} and also the regions related to the exponential holding BM and absorbed BM. This is natural to expect since in the killed BM setting there is an equivalence between the origin and the cemetery~$\Delta$.
	\end{remark}
	
	\section{Simulations}\label{s4}
	In this section we provide computational simulations of the boundary random walk's path. The Python script which generated them is available  at the GitHub repository \cite{BRW}. Sample outputs of this code  corresponding to each possible limit of the boundary random walk are showed in Figures \ref{fig:elastic}, \ref{fig:sticky},  \ref{fig:exponential_holding},  \ref{fig:reflected}, \ref{fig:absorbed}, \ref{fig:killed}  and  \ref{fig:mixed}.
	
	\newcommand{\constante}{0.45}

	\begin{figure}[H]
		\centering
		\begin{subfigure}[t]{\constante\textwidth}
			\includegraphics[width=\linewidth]{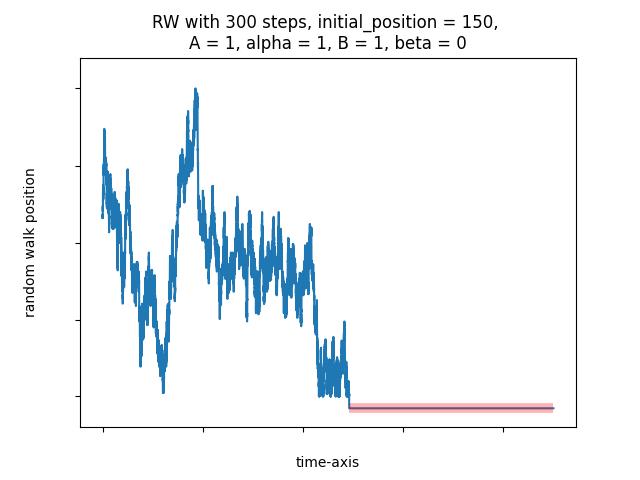}
			\caption{\textbf{Elastic BM.} Here we took $\alpha=1$ and  $\beta = 0$. Under this choice of parameters, the random walk  approximates  the elastic BM, which is a mixture of being reflected and killed. The cemetery state is highlighted in red.}
			\label{fig:elastic}
			
		\end{subfigure}
		\hfill
		\begin{subfigure}[t]{\constante\textwidth}
			\includegraphics[width=\linewidth]{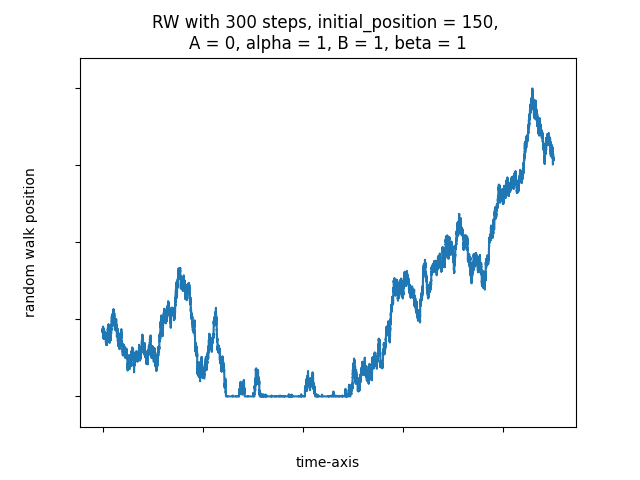}
			\caption{\textbf{Sticky BM.} Here we took $\beta = 1$ and $A = 0$. Under this choice of parameters, the random walk  approximates  the sticky Brownian motion. As we can see above, once the walk touches zero, it takes a certain time to scape from the stickiness of the origin.}
			\label{fig:sticky}
		\end{subfigure}
		\caption{Convergence to Elastic  and Sticky BM}
		\label{fig:figura_principal6}
	\end{figure}

	\begin{figure}[H]
		\centering
		\begin{subfigure}[t]{\constante\textwidth}
			\includegraphics[width=\linewidth]{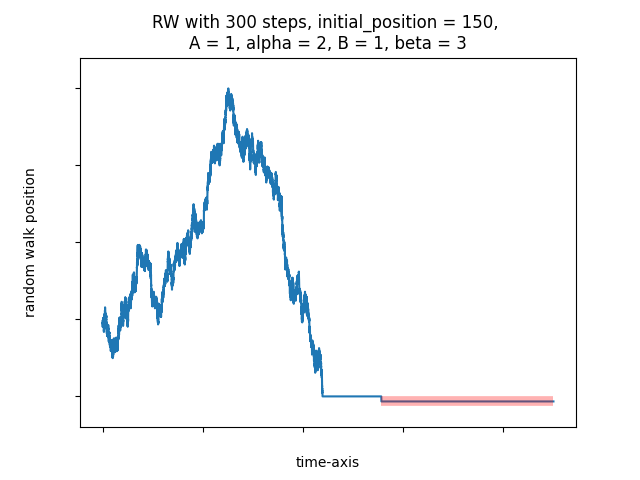}
			\caption{\textbf{Exponential Holding BM.} Here we took $\alpha = 2$ and  $\beta = 3$. Under this choice of parameters, the random walk  approximates  the exponential holding Brownian motion, where the walk after reaching zero, stays there for an exponential time and then jumps to the cemetery. The cemetery state is highlighted in red.}
			\label{fig:exponential_holding}
		\end{subfigure}
		\hfill
		\begin{subfigure}[t]{\constante\textwidth}
			\includegraphics[width=\linewidth]{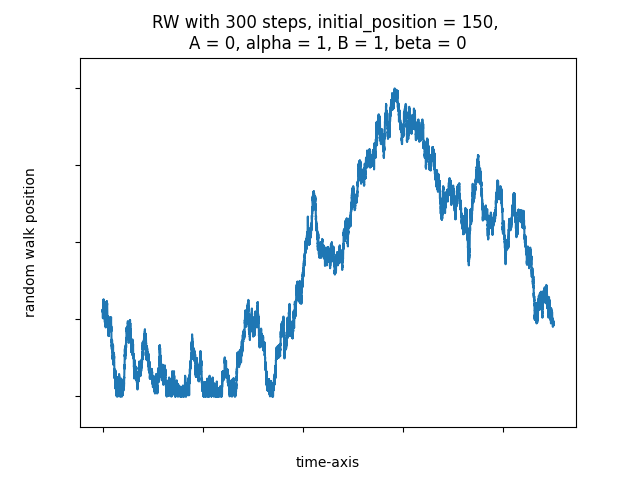}
			\caption{\textbf{Reflected BM.}  Here we took $\beta = A = 0$. Under this choice of parameters, the random walk  approximates  the  reflected Brownian motion. Note that the absence of stickiness (compare it  with Figure~\ref{fig:sticky}).}
			\label{fig:reflected}
		\end{subfigure}
		\caption{Convergence to Exponential Holding and Reflected Brownian Motion}
		\label{fig:figura_principal7}
	\end{figure}

	\begin{figure}[H]
		\centering
		\begin{subfigure}[t]{\constante\textwidth}
			\includegraphics[width=\linewidth]{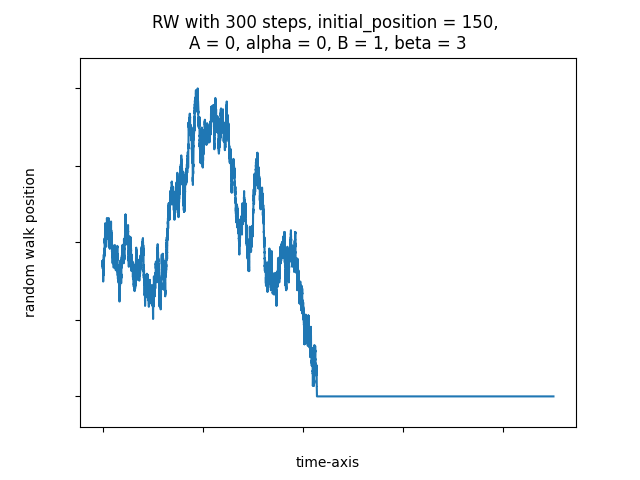}
			\caption{\textbf{Absorbed BM.} Here $A = 0$ and $\beta = 3$. Under this choice of parameters, the random walk approximates the absorbed Brownian motion. Note that once the walk reaches the origin, it stays there forever.}
			\label{fig:absorbed}
		\end{subfigure}
		\hfill
		\begin{subfigure}[t]{\constante\textwidth}
			\includegraphics[width=\linewidth]{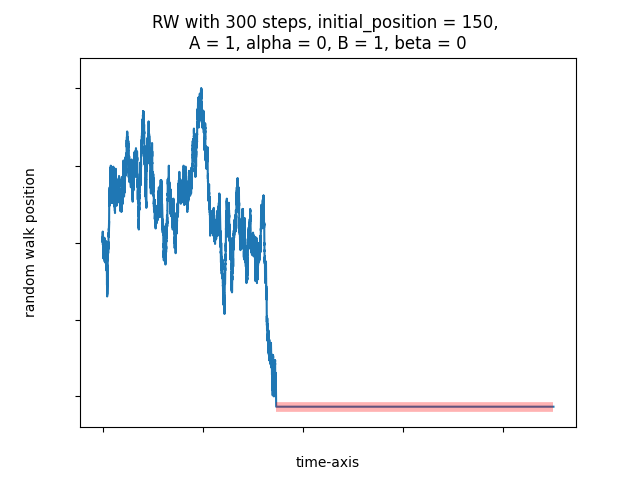}
			\caption{\textbf{Killed BM.} Here we took
				$\alpha = \beta = 0$ and $A = B = 1$.  Under this choice of parameters, the random walk  approximates  the killed BM, where the process, once reaching zero, immediately jumps to the cemetery state, which is highlighted in red.}
			\label{fig:killed}
		\end{subfigure}
		\caption{Convergence to Absorbed and Killed Brownian Motion}
		\label{fig:figura_principal8}
	\end{figure}

	\begin{figure}[H]
		\includegraphics[width=\constante\linewidth]{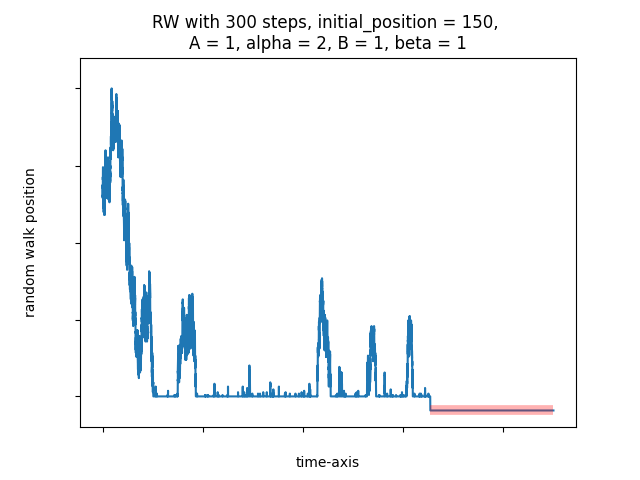}
		\caption{\textbf{Mixed Brownian Motion.} Here  $\alpha = 2$ and  $\beta = 1$.  Under this choice of parameters, the random walk  approximates  the mixed Brownian motion. Note that the process has some stickiness at zero, and a certain random time jumps to the cemetery state, which is highlighted in red.}
		\label{fig:mixed}
	\end{figure}

	\section*{Funding}
	D.E.\ was supported by the National Council for Scientific and Technological Development - CNPq via a Bolsa de Produtividade  318182/2025-4. D.E.\ moreover acknowledge support by the Serrapilheira Institute (Grant Number Serra-R-2011-37582). D.E.\ and T.F.\ were  partially supported by FAPESB (EDITAL FAPESB  012/2022 - UNIVERSAL - APP0044/2023). FAPESB is the Bahia Research Foundation. D.E. and
	T.F. acknowledge support by the CNPq via a Universal Grant (Grant Number 401314/2025-
	1). D.E, T.F.\ and M.J.\  acknowledge support by the National Council for Scientific and Technological Development - CNPq via a Universal Grants (Grant Numbers 406001/2021-9 and 401314/2025-1). T.F.\ was supported by the National Council for Scientific and Technological Development - CNPq via a Bolsa de Produtividade number 311894/2021-6. E.P.\ thanks FAPESB (Funda\c c\~ao de Amparo \`a Pesquisa do Estado da Bahia) for the support via a PhD scholarship.

	\appendix
			\section{Local Lipschitz property of semigroups}\label{appendixA}
This appendix provides a verification of the local Lipschitz property of the semigroups of the absorbed Brownian motion, killed Brownian motion,
	reflected Brownian motion, sticky Brownian motion, elastic Brownian motion, exponential holding Brownian motion and mixed Brownian motion. This is needed to verify condition \ref{(G4)}, which is an ingredient in the proof of Theorem~\ref{thm2.3}.
	 Namely, it  says that $\{T(t) : t \ge 0\}$ is locally Lipschitz in the sense that, for each $t>0$ and each compact set $K\subseteq \msf{S}$, there exists a constant $M = M(t,K)>0$ such that, for all $f\in \mc{C}_0(\msf S_\Delta)$,
	\begin{equation*}
		\vert  \msf{T}(t)f (u) - \msf{T}(t)f (\tilde{u}) \vert \; \leq \; M\cdot \Vert f\Vert_\infty \cdot  d_{\msf{S}}(u,\tilde{u})\,,\quad \forall\, u,\tilde{u}\in K\,.
	\end{equation*}
	For the processes considered in this paper,  we have $\msf{S} = [0,\infty)$, except  in the case of the killed BM, where $\msf{S} = (0,\infty)$. In all cases, the distance on $\msf{S}$ is $d(u,\tilde{u}) = |u-\tilde{u}|$. Let us postpone the analysis of  the killed BM to the end. The semigroups of absorbed BM, 
	reflected BM, sticky BM, elastic BM, exponential holding BM and mixed BM  all have  transition probability densities of the form 
	\begin{equation*}
		p_t(u, dv) \;=\; q_t(u,v) dv + r_t(u)\delta_0(dv) + s_t(u)\delta_\Delta(dv)\,.
	\end{equation*}
	Since the test functions are assumed to satisfy $f(\Delta)=0$, we can estimate 
	\begin{align*}
		\vert  \msf{T}(t)f (u) - \msf{T}(t)f (\tilde{u}) \vert \; \leq \; &  \int_{[0,\infty)} |q_t(u,v)- q_t(\tilde{u}, v)| \cdot |f(v)| \,dv  +|r_t(u)- r_t(\tilde{u})| \cdot |f(0)|  \\
		\;\leq\; & \Vert f\Vert_\infty \Big\{|r_t(u)- r_t(\tilde{u})|+  \int_{[0,\infty)} |q_t(u,v)- q_t(\tilde{u}, v)| \,dv   \Big\} \,.
	\end{align*}
	Assume for a moment that we have checked the following two conditions:
	\begin{enumerate}[label=(A\arabic*)]
		\item \label{item:lemLIP1} The function $r_t(u)$ is $C^1$ in $u\in [0,\infty)$.
		\item \label{item:lemLIP2} The function $q_t(u,v)$ is $C^1$ in $u\in [0,\infty)$ and, for all $K\subset \msf{S}$ compact, it holds that  
		\begin{equation*}
			\int_{[0,\infty)} \Big[\sup_{u\in K}|\p_u q_t(u,v)|\Big]  dv \;<\;\infty\,.
		\end{equation*}
	\end{enumerate}
	Then, under \ref{item:lemLIP1} and \ref{item:lemLIP2}, we have 
	\begin{align*}
		\vert  \msf{T}(t)f (u) - \msf{T}(t)f (\tilde{u}) \vert 
		\;\leq\; & \Vert f\Vert_\infty \cdot \bigg\{\sup_{u\in K} \vert \partial_u r_t(u,v)\vert +  \int_{[0,\infty)}\Big[\sup_{u\in K} |\p_u q_t(u,v)|\Big] \,dv   \bigg\} \cdot \vert u-\tilde{u}\vert \\
		\; = \; &\Vert f\Vert_\infty \cdot M(t, K) \cdot \vert u-\tilde{u}\vert\,,
	\end{align*}
	assuring the local Lipschitz property. So our task is to check \ref{item:lemLIP1} and \ref{item:lemLIP2} for each of the seven semigroups we dealt with in this paper. Except for exponential holding BM, all formulas below are adapted from \cite[Appendix~1]{handbook}. Note that there is missing a factor of $1/2$ everywhere in comparison with the formulas of \cite{handbook} because there  the authors take  twice the Lebesgue measure as reference measure.  Denote by 
	\begin{equation*}
		\Gamma(z) \;=\; \frac{2}{\sqrt{\pi}}\int_z^\infty e^{-u^2}\dd u \qquad \text{ and } \qquad G_t(u,v) \;=\; \frac{e^{-\frac{(u-v)^2}{2t}}}{\sqrt{2\pi t}} 
	\end{equation*}
	the error function and the Gaussian kernel, respectively. 
	Note that these two functions satisfy the following properties.
	\begin{enumerate}[label=(\roman*)]
		\item \label{(i)} $\Gamma$ and $G_t(\cdot,v)$ are smooth for all times $t>0$ and all $v$. 
		\item \label{(ii)} There exists a $c>0$ such that for all $z>0$ one has that $\Gamma(z)\leq \tfrac{c}{z} e^{-z^2}$. This is a known Gaussian tail bound.
		\item \label{(iii)} A direct computation yields $\Gamma'(z)=-\tfrac{2}{\sqrt{\pi}}e^{-z^2}$.
		\item \label{(iv)} For all compact sets $K\subset \msf{S}$ and all times $t>0$ the functions
		\begin{equation*}
			v\mapsto \sup_{u\in K} |\partial_u G_t(u,v)|\qquad\text{and}\qquad v\mapsto \sup_{u\in K} |\partial_u G_t(u,-v)|
		\end{equation*}
		are integrable over $v\in[0,\infty)$. This follows from a direct computation.
	\end{enumerate}

	\textbf{Absorbed BM} \cite[Page 121]{handbook}. In this case,  
	\begin{align*}
		{q}_t(u,v) & \;=\; G_t(u,v) -  G_t(u,- v)
		\qquad \text{ for }u,v\in [0,\infty),\\
		{r}_t(u) & \;=\; 1-  \int_0^\infty q_t(u,v)dv, \quad \text{ for } u\in [0,\infty),\\
		{s}_t(u) & \;=\; 0.
	\end{align*}
	
	\textbf{Reflected BM} \cite[Page 120]{handbook}. In this case, 
	\begin{align*}
		{q}_t(u,v) & \;=\; G_t(u,v) +  G_t(u,- v),
		\qquad \text{ for }u,v\in [0,\infty),\\
		{r}_t(u) & \;=\; 0,\\
		{s}_t(u) & \;=\; 0.
	\end{align*}

	\textbf{Elastic BM} \cite[Page 125]{handbook}. In this case, denoting $\gamma := \frac{c_1}{c_2}$, where $c_1, c_2$ are as in \citealp[Theorem 2.10]{EFJP}, we have that
	\begin{align*}
		{q}_t(u,v) & \;=\; G_t(u,v) +  G_t(u,- v) - \gamma e^{\Big(\gamma(u + v) + \pfrac{\gamma^2t}{2}\Big)}\Gamma\Bigl(\pfrac{u + v + \gamma t}{\sqrt{2t}}\Bigr),
		\qquad \text{ for }u,v\in [0,\infty),\\
		{r}_t(u) & \;=\; 0, \\
		{s}_t(u) & \;=\;  1-  \int_0^\infty q_t(u,v)dv, \quad \text{ for } u\in [0,\infty).
	\end{align*}
	Note that the transition kernel in \cite{handbook} slightly differs from ours. This is due to the fact that the transition density $p$ in \cite{handbook} is sub-Markovian. Our choice of $s$ turns the kernel into a Markovian one.

	\textbf{Sticky BM} \cite[Pages 123-124]{handbook}. In this case,
	\begin{align*}
		{q}_t(u,v) & \;=\; G_t(u,v) -  G_t(u,- v) + \frac{1}{\gamma}e^{\Big(\pfrac{2u+ 2v}{\gamma} + \pfrac{2t}{\gamma^2}\Big)}\Gamma\Bigl(\pfrac{u+ v}{\sqrt{2t}} + \pfrac{\sqrt{2t}}{\gamma}\Bigr)\,,
		\;\text{ for }u,v\in [0,\infty),\\
		{r}_t(u) & \;=\; 1-  \int_0^\infty q_t(u,v)dv, \quad \text{ for } u\in [0,\infty),\\
		{s}_t(u) & \;=\; 0.
	\end{align*}
	Note that in \cite{handbook} the sticky Brownian motion is defined on the full line. Nevertheless, identifying $-x$ with $x$ one can deduce the corresponding formula on the half-line.

	\textbf{Exponential Holding BM}. In this case, the transition probability density is 
	\begin{align}
		{q}_t(u,v) & \;=\; G_t(u,v) -  G_t(u,- v),
		\qquad \text{ for }u,v\in [0,\infty), \notag \\
		{r}_t(u) & \;=\; 
		\begin{cases}
			\int_{0}^t f_{\tau_0}(s,u) e^{-\lambda (t-s)}ds &  \text{ for } u\in(0,\infty),\vspace{2pt}\\
			e^{-\lambda t} & \text{ for } u=0,
		\end{cases} \label{eq_exp}\\
		{s}_t(u) & \;=\; 
		\begin{cases}
			\int_{0}^t f_{\tau_0}(s,u) \big(1- e^{-\lambda (t-s)}\big)ds, & \text{ for } u\in (0,\infty),\\
			1- e^{-\lambda t} & \text{ for } u=0,\notag
		\end{cases}
	\end{align}
	where  
	$f_{\tau_0}(s,u) = \frac{u}{\sqrt{2\pi s^3}}e^{-\frac{u^2}{2s}}\one_{[0,\infty)}(s)$ is the density of $\tau_0$, the hitting time of zero for the absorbed BM.
	We were not able to locate the above case in the literature. However one can deduce it as follows. The exponential holding Brownian motion is defined as a Brownian motion that upon hitting zero cannot transition to the half-line anymore (so it gets absorbed), and then jumps from zero to the cemetery after an exponentially distributed amount of time. The formula for $q_t$ coincides with the corresponding expression for the absorbed Brownian motion. The first line in the definition of $r_t$ integrates the probability density of the process to go from $u$ to zero at time $s$ and then staying there for the remaining $t-s$ time units. The second line is the probability of an exponentially distributed random variable to be bigger than $t$, which means that the process starts in zero and stays there. An analogous reasoning can be applied to deduce the formula for $s_t$. \bigskip

	\textbf{Mixed BM} \cite[Page 125]{handbook}. In this case, denoting $\gamma = \frac{c_1}{c_2}$,   $c = \frac{c_3}{2c_2}$ and $\Upsilon = \sqrt{1 - \gamma c}$,  where $c_1, c_2, c_3$ are as in \citealp[Theorem 2.10]{EFJP}, we have that
	\begin{align*}
		{q}_t(u,v)  \;=\; &  G_t(u,v) -  G_t(u,- v)  + \pfrac{\Upsilon + 1}{2 c \Upsilon}e^{\big(\pfrac{(u + v)2\gamma}{1 - \Upsilon} + \pfrac{(\Upsilon + 1)^2 t}{8 c^2}\big)}\Gamma\bigg( \pfrac{u + v}{\sqrt{2t}} + \pfrac{(1 + \Upsilon)\sqrt{t}}{c\sqrt{2}}\bigg) \\
		& + \pfrac{\Upsilon - 1}{2 c \Upsilon} e^{\big( \pfrac{(u + v)2\gamma}{1 + \Upsilon} + \pfrac{(\Upsilon - 1)^2 t}{8 c^2}\big)} \Gamma\bigg(\pfrac{u + v}{\sqrt{2t}} + \pfrac{(1 - \Upsilon)\sqrt{t}}{2c\sqrt{2}}\bigg), 	\qquad \text{ for }u,v\in [0,\infty),\\
		{r}_t(u)  \;=\;  &c q_t(u,0), \quad \text{ for } u\in [0,\infty),\\
		{s}_t(u)    \;=\;  & 1-  r_t(u) - \int_0^\infty q_t(u,v)dv, \quad \text{ for } u\in [0,\infty).
	\end{align*}

	For the absorbed, and reflected Brownian motion, conditions \ref{item:lemLIP1} and \ref{item:lemLIP2} are consequences of items \ref{(i)} and \ref{(iv)} above together with dominated convergence. For the elastic, sticky and mixed BM one additionally needs to use items \ref{(ii)} and \ref{(iii)} above to control the multiplications of $\Gamma$ by the various exponentials. Finally, for the exponential holding Brownian motion, the condition~\ref{item:lemLIP2} follows by the same arguments as above. Regarding the condition \ref{item:lemLIP1} we will show by a direct calculation that the derivative of $r_t$ is uniformly bounded for all $u>0$ and then we establish the Lipschitz property in $u=0$ by a coupling argument. We start with the case $u>0$.
	Note that we need to analyse the derivative of
	\begin{equation*}
		\int_0^t \frac{u}{s^{3/2}} e^{-u^2/2s} e^{\lambda s}\, ds\,.
	\end{equation*} 
	Using the substitution $v=\tfrac{u^2}{2s}$, we see that the above integral equals
	\begin{equation*}
		\sqrt{2}\int_{\tfrac{u^2}{2t}}^\infty v^{-1/2} e^{-v} e^{\lambda \tfrac{u^2}{2v}}\, dv\,.
	\end{equation*}
	To continue, define for $x>0$
	\begin{equation*}
		h(x,u)\;=\; \int_x^\infty v^{-1/2} e^{-v} e^{\lambda \tfrac{u^2}{2v}}\, dv\,,
	\end{equation*}
	and note that
	\begin{equation*}
		\partial_u \big(v^{-1/2} e^{-v} e^{\lambda \tfrac{u^2}{2v}}\big)\;=\; v^{-1/2} e^{-v} \tfrac{u}{v}e^{\lambda \tfrac{u^2}{2v}}\,.
	\end{equation*}
	Hence, dominated convergence implies that we can interchange integration and differentiation with respect to $u$ in the integral above.
	Using the chain rule we thus conclude that
	\begin{equation*}
		\partial_u h(\tfrac{u^2}{2t}, u)
		\;=\; \frac{\sqrt{2}}{\sqrt{t}} e^{-\tfrac{u}{2t}}e^{\lambda t}
		+ \int_{\tfrac{u^2}{2t}}^\infty \frac{u}{v^{3/2}} e^{-\lambda v}e^{\lambda \tfrac{u^2}{2v}}\, dv\,.
	\end{equation*}
	Estimating $e^{\lambda \tfrac{u^2}{2v}}\leq e^{\lambda t}$ on the domain of integration we see that, for all $u>0$:
	\begin{equation*}
		\int_1^\infty \frac{u}{v^{3/2}} e^{-\lambda v}e^{\lambda \tfrac{u^2}{2v}}\, dv\;\leq\; Cu\,,
	\end{equation*}
	for a constant $C>0$ that does not depend on $u$.
	Finally (assuming $u^2\leq 2t$ otherwise there is nothing to prove) and estimating $e^{-\lambda v}$ from above by one, we can estimate
	\begin{equation*}
		\int_{\tfrac{u^2}{2t}}^1 \frac{u}{v^{3/2}} e^{-\lambda v}e^{\lambda \tfrac{u^2}{2v}}\, dv
		\;\leq\; 	e^{\lambda t}\int_{\tfrac{u^2}{2t}}^1\frac{u}{v^{3/2}}\, dv
		\;=\; e^{\lambda t}\frac32 (\sqrt{2t}-u)\,.
	\end{equation*}
	This concludes the prove that $\partial_u r_t(u)$ is bounded uniformly in $u>0$. 
	
	We turn to the case $u=0$.
	Denote by $B^u_t$ the exponential holding BM  starting from $u\geq 0$. Fix $u>0$. First, we construct a  coupling of $B^0_t$ and $B^u_t$. The exponential holding BM is the process that behaves as the absorbed BM until it hits zero. After that, it stays at zero for an exponential time of parameter $\lambda>0$ and then jumps to the cemetery $\Delta$. Let $\xi_1, \xi_2,\ldots$ be the arrivals of a Poisson  process of parameter $\lambda$ and let $W_t^u$ be a standard BM on the real line starting from $u$, independent of that Poisson process, and let $\tau_0$ be its hitting time of zero. Let $N = \inf\{k\geq 1: \xi_k \geq \tau_0\}$ and define
	\begin{align*}
		B^0_t \;=\; \begin{cases}
			0, & \text{ if } t< \xi_1,\\
			\Delta, & \text{ if } t\geq \xi_1,
		\end{cases}
		\qquad 	\text{and} \qquad
		B^u_t \;=\; \begin{cases}
			W_t^u, & \text{ if } t< \tau_0,\\
			0, & \text{ if } \tau_0 \leq t\leq \xi_N,\\
			\Delta, & \text{ if }  t\geq  \xi_N.
		\end{cases}   
	\end{align*} 
	In words, $B^0_t$ and $B^u_t$ use the same Poisson process to jump to the cemetery, see the Figure~\ref{fig:exp_coupling} below for an illustration.
	Note that by the independence of the Brownian motion and the Poisson process and the memoryless property of the latter, conditionally on $\tau_0<\infty$, the random variable $\xi_N-\tau_0$ is exponentially distributed with parameter $\lambda$ and is independent of the Brownian path up to time $\tau_0$. Hence, $B^u$ evolves as Brownian motion until it first hits zero, remains at zero for an independent exponential time of parameter $\lambda$, and afterwards jumps to the cemetery state. This is precisely the exponential holding Brownian motion. See also \cite[Sections 3 and 5]{Kostrykin} for the corresponding pathwise construction.
	
	\newcommand{\browniano}[4]{
		\draw[#4, thick] (0,0.85)
		\foreach \x in {1,...,#1}
		{   -- ++(#2,-rand*#3-.006)
		};
	}
	\begin{figure}[!htb]
		\centering
		\begin{tikzpicture}[scale=1];
			\pgfmathsetseed{5}
			\browniano{112}{0.02}{0.13}{red}{}
			\draw[fill=black] (0,0.85) circle (1.5pt) node[left]{$u$};
			
			\draw (0,0) node[anchor = north east]{$0$};
			
			\draw[->, thick] (0,-1.25)--(0,1.75);
			\draw[->] (-1,0)--(8,0) node[anchor=north]{time};
			\draw[very thick, red] (2.22,0) -- (3.5,0);
			\draw (2.22,0) node[below]{$\tau_0$};
			\draw[very thick, red] (3.5,-1) -- (7.5,-1);
			\draw[dashed] (1.5,0)-- (1.5,-1);
			\draw[dashed] (3.5,0)-- (3.5,-1);
			\draw[dashed] (0,-1)-- (3.5,-1);
			\draw [very thick, blue] (0,0) -- (1.5,0);
			\draw[red, fill = white] (3.5,0) circle (1.5pt) node[below left, black]{$\xi_2$};
			\draw[black, fill = black] (6,0) circle (1.5pt) node[below, black]{$\xi_3$};
			\draw[blue, fill = white] (1.5,0) circle (1.5pt) node[below left, black]{$\xi_1$};
			\draw [very thick, blue] (1.5,-1) -- (3.5,-1);
			\draw[red, fill = red, thick] (3.5,-1) circle (1.5pt);
			\draw[blue, fill = blue, thick] (1.5,-1) circle (1.5pt);
			\draw[fill=black] (0,-1) circle (1.5pt) node[left]{$\Delta$};	
		\end{tikzpicture}
		\caption{Illustration of the coupling between $\{B^u_t:t\geq 0\}$ and $\{B^0_t:t\geq 0\}$. In red, the trajectory of $\{B^u_t:t\geq 0\}$ and in blue, the trajectory of $\{B^0_t:t\geq 0\}$.  In the picture above, $N=2$. For $u$ small, we typically  have that $\tau_0< \xi_1$ and, in this case, $B^u_t = B^0_t$ for $t>\tau_0$.}\label{fig:exp_coupling}
	\end{figure}
	Denote $A_t = [\tau_0<t]\cap [\tau_0< \xi_1]$. Then,
	\begin{align*}
		\big\vert \bb E\big[f(B^u_t)\big] - \bb E\big[f(B^0_t)\big]\big\vert & \;\leq\;   \bb E\big[ \big\vert f(B^u_t) - f(B^0_t) \big\vert \cdot\one_{A_t}\big] +
		\bb E\big[ \big\vert f(B^u_t) - f(B^0_t) \big\vert \cdot \one_{A_t^\complement}\big]\\
		& \; =\;   \bb E\big[ \big\vert f(B^u_t) - f(B^0_t) \big\vert \cdot \one_{A_t^\complement}\big] \\
		& \;\leq\;   2 \Vert f \Vert_\infty \cdot \bb P (A_t^\complement)\\
		& \;\leq\;   2 \Vert f \Vert_\infty \cdot \big\{\bb P (\tau_0 \geq t) +\bb P (\tau_0 \geq \xi_1)\big\}\,. 
	\end{align*}
	By the Reflection Principle, we know that $
	\bb P (\tau_0 \geq t)= \int_{-u}^u \frac{e^{-\frac{x^2}{2t}}}{\sqrt{2\pi t}}dx \;\leq\; \frac{2u}{\sqrt{2\pi t}}$.
	Since $\tau_0$ and $\xi_1$ are independent, by the Substitution Principle (see \citealp[Example~5.1.5]{Durrett}), we have that
	\begin{align*}
		\bb P (\tau_0 \geq \xi_1)  & \; = \; \bb E \Big[ \bb E \big[\one_{\tau_0 \geq \xi_1} \mid \xi_1 \big]\Big]\; \leq  \;  \bb E\Big[ \frac{2u}{\sqrt{2\pi \xi_1}}\Big] \;=\; u \cdot \bb E\Big[ \frac{2}{\sqrt{2\pi \xi_1}}\Big]\,.
	\end{align*}
	Since $\xi_1\sim \mathrm{exp}(\lambda)$, then $\bb E[1/\sqrt{\xi_1}]<\infty$, and we conclude that the semigroup of the exponential holding BM is Lipschitz at zero.
	\bigskip
	
	\textbf{Killed BM} \cite[Pages 120 - 121]{handbook}. In this  case, $\msf{S} = (0,\infty)$, where $0$ and $\Delta$ are identified.  Its semigroup is given by
	\begin{align*}
		{q}_t(u,v) & \;=\; G_t(u,v) -  G_t(u,- v)
		\qquad \text{ for }u,v\in [0,\infty),\\
		{r}_t(u) & \;=\; 1-  \int_0^\infty q_t(u,v)dv, \quad \text{ for } u\in [0,\infty),
	\end{align*}
	where the functions $f$ in its domain must satisfy $f(0)=0$, and \ref{item:lemLIP2} is enough to guarantee the local Lipschitz property, which is a consequences of items \ref{(i)} and \ref{(iv)} together with dominated convergence. Note that the killed Brownian motion essentially coincides with the absorbed Brownian motion. The difference however is that $0$ is identified with the cemetery.

	
	\bibliographystyle{plainnat} 
	\bibliography{bibliografia}       
	
	
\end{document}